\RequirePackage{fix-cm}
\documentclass{article}       
\usepackage{graphicx}
\RequirePackage{amsfonts,amssymb,stmaryrd}
\RequirePackage{mathtools, xcolor}
\usepackage{breqn, siunitx, algpseudocode}
\DeclareSIUnit\year{yr}
\RequirePackage{bm, graphicx, xspace}
\RequirePackage{soul}
\usepackage{commands}
\usepackage{ulem}
\usepackage{caption}
\usepackage{subcaption}
 \usepackage{hyperref}
 \hypersetup{
    colorlinks=true,       
    hidelinks=false,
    linkcolor=red,          
    citecolor=green,        
    filecolor=cyan,         
    urlcolor=magenta        
}
 
\graphicspath{{figures/}}
\begin{document}

\title{Numerical modelling of convection-driven cooling, deformation and fracturing of thermo-poroelastic media}
%



\author{Ivar Stefansson \and Eirik Keilegavlen \and Sæunn Halldórsdóttir \and Inga Berre 
}



\date{}

\maketitle

\begin{abstract}
Convection-driven cooling in porous media influences thermo-poro-mechanical stresses, thereby causing deformation.
These processes are strongly influenced by the presence of fractures, which dominate flow and heat transfer. At the same time, the fractures deform and propagate in response to changes in the stress state. Mathematically, the model governing the physics is tightly coupled and must account for the strong discontinuities introduced by the fractures. Over the last decade, and motivated by a number of porous media applications, research into such coupled models has advanced modelling of processes in porous media substantially. 

Building on this effort, this work presents a novel model that couples flow, heat transfer, deformation, and propagation of fractures with flow, heat transfer, and thermo-poroelasticity in the matrix. The model is based on explicit representation of fractures in the porous medium, and discretised using multi-point finite volume methods. Frictional contact and non-penetration conditions for the fractures are handled through active set methods,
while a propagation criterion based on stress intensity factors governs fracture extension. Considering both forced and natural convection processes, the numerical results show the intricate nature of thermo-poromechanical fracture deformation and propagation. 

\end{abstract}

%
%

\section{Introduction}
For a porous medium, possibly containing fractures, the interplay between flow, thermal transport, and deformation can be strong. 
In particular, cooling of the medium induces thermal stress that can lead to deformation and fracturing.
Furthermore, fractures deform and propagate as a result of the coupled dynamics. The result is coupled thermo-hydro-mechanical (THM) processes in the intact porous medium, interacting with flow and thermal transport in fractures as well as fracture deformation and propagation. Such coupled process-structure interaction is characteristic for a wide range of natural and engineered processes in natural and manufactured materials. For example, the structural and functional performance of concrete structures, like dams, bridges, nuclear and liquefied natural gas containers, and cement sheaths of subsurface well bore constructions, are affected by the time evolution of their properties under variable THM loads \cite{baroth_uncertainty_2019,bois2012use,kogbara2013review,bouhjiti_statistical_2018,bouhjiti_statistical_2018-1,lin_integrity_2020}. In the subsurface, THM processes interact with deformation and propagation of fractures in fluid injection operations \cite{siratovich2015saturated,gao_three-dimensional_2020,wu_investigation_2020}. The coupled dynamics is also hypothesised to be crucial in heat transfer from the deep roots of geothermal systems by deepening natural convection through evolving fractures \cite{lister1974penetration,bodvarsson1982terrestrial,bjornsson_penetration_1982,bjornsson_heatmass1987}. Common to all these applications is that tight coupling in the dynamics limits the knowledge which can be gained from analysis of individual processes and mechanisms in isolation. This motivates development of simulation models that acknowledge the coupled nature of the physics. 

Since its foundation by Biot \cite{biot1941general}, the theory of poroelasticity has successfully been applied to model coupled hydro-mechanical processes.
The extention to thermo-poroelasticity \cite{coussy2004poromechanics} is also widely applied, including in geomechanics \cite{selvadurai2017thermo}. More recently, models accounting for discontinuities in the form of fractures in poroelastic and  thermo-poroelastic media have been developed.
Typically, the development has focused either on deformation of preexisting fractures or the mechanical fracturing of the materials.
The models can be distinguished based on whether fractures are represented explicitly as discrete objects embedded in the porous medium, or represented as part of the porous medium itself.
The latter incorporate the effect of the extent to which the material is fractured by use of smeared or distributed representations. Such models include phase-field and damage approaches for fracture \cite{de2016gradient} and continuum and multi-continuum approaches for flow models \cite{berre2019flow}.

Approaches based on explicit representation of the fractures can further be distinguished by how the fractures are represented in discretisation, 
specifically on whether a conforming or non-conforming representation of the fractures is used in the grid \cite{berre2019flow}. 
Non-conforming methods represent the fracture through an enriched representation. For poromechanics, combinations of the embedded discrete fracture method, extended finite element methods and/or embedded finite element methods have been applied \cite{ren2016fully,cusini2020simulation,giovanardi2017unfitted}.
Such non-conforming approaches have also been extended to include tensile fracture propagation based on extended finite element \cite{khoei2014mesh} and embedded discrete fracture methods \cite{deb2020extended}. 
Conforming methods use a representation where the fractures coincide with matrix faces. Considering fractures that have a negligible aperture compared to the modelled domain, this representation can be combined with an approach where fractures are modelled as lower-dimensional structures 
\cite{Martin2005,Karimi-Fard2003efficient} and discretised with elements of zero thickness \cite{kiraly1988mixedfem,Boon2018,flemisch2018benchmarks,berre2020verification}. 
For poroelastic media with fractures, this allows for the application of standard finite element \cite{salimzadeh2017three}, finite volume \cite{ucar2018finite,berge2020hmdfm} and combined finite element/finite volume schemes \cite{garipov2016discrete,settgast2017fully,garipov2019discrete}, more recently also including fracture contact mechanics  \cite{garipov2016discrete,gallyamov2018discrete,franceschini2020algebraically}
and tensile fracture propagation \cite{settgast2017fully,salimzadeh2017three}.
Thermal effects on fracture deformation and propagation in poroelastic media are less studied, although some recent studies model deformation of existing fractures in thermo-poroelastic media \cite{salimzadeh2018thm,stefansson2020thm,garipov2019discrete}.

Motivated by the development of increasingly sophisticated THM models for fractured porous media, our goal in the present paper is to extend numerical modeling of THM to also incorporate fracture propagation, and thereby contribute to bridge the gap between fracture mechanics models and coupled THM models for porous media. 
Specifically, we consider mathematical and numerical modelling of fracture deformation and propagation resulting from coupled THM-processes. 
Our focus is on convection-driven cooling in the subsurface, where forced or natural fluid convection induces thermo-poromechanical stress changes leading to fracture deformation and propagation.
The dynamics is characterised by tight coupling between physical processes and strong interaction between the physical processes and the (evolving) geometry of the fracture network. 
Accordingly, our model and simulation approach is designed to faithfully represent these couplings, including fracture deformation and propagation.

The fractured thermo-poroelastic medium is represented using a discrete fracture-matrix model, where fractures are represented as lower-dimensional discontinuities in an otherwise continuous thermo-poroelastic medium.
Deformation of existing fractures is modelled through contact mechanics relations based on a Coulomb friction criterion for slip along the fractures and a non-penetration condition \cite{hueber2008paper,berge2020hmdfm}.
This is combined with a simple criterion for fracture propagation based on the mode I stress intensity factor, which we compute directly from the displacement jump in the vicinity of the fracture tip using a variant \cite{nejati2015} of the displacement correlation method \cite{chan1970displacementcorrelation}. To adjust the grid to an arbitrary fracture propagation path is highly technical \cite{paluszny2011,Hau2020}, and we instead
make the assumption that fractures propagate along existing faces in the matrix grid.
This constrains the numerical representation of an evolving fracture and makes it difficult to preserve reasonable fracture geometries for general propagation scenarios, in particular for three-dimensional problems.
We therefore further limit ourselves to tensile fracturing, where the possible propagation path is easy to predict and the grid can be constructed to accommodate the propagation.

We discretise the model using a control volume framework for fracture contact mechanics in thermo-poroelastic media \cite{stefansson2020thm}. 
The control volume approach builds on a combination of the multi-point stress approximation method for Biot poroelasticity \cite{nordbotten2016biot,keilegavlen2017finite} with the multi-point flux approximation method for flow \cite{Aavatsmark2002}. 
This combination is previously applied for numerical modelling of fractured poroelastic media \cite{ucar2018finite} with a simplified model for deformation along fractures. 
The fracture contact mechanics builds on work by Berge et al. \cite{berge2020hmdfm}, who formulated the contact conditions on the fracture using Lagrange multipliers representing the contact tractions \cite{wohlmuth2011}. 
Using this approach, the variational inequality representing the contact problem can be rewritten using complementary functions, and the resulting system of equations solved by a semi-smooth Newton method \cite{hueber2008paper,berge2020hmdfm}.
Our model is implemented in the open-source simulator PorePy \cite{keilegavlen2020porepy}, which is designed for multiphysics problems in fractured porous media. 

We assess the reliability of our simulation tool by tests that probe the approximations of both the onset of fracturing and the speed of fracture propagation.
We then present two application-related simulations that both involve fracture propagation driven by convective cooling. 
The cases include respectively forced convection during production of geothermal energy and natural convection in vertical fractures in the presence of high thermal gradients.
Taken together, the results show the importance of  developing simulation tools that can accurately represent the tight couplings in THM processes, and also deal with deformation and propagation of fractures.

The paper is structured as follows. Section \ref{sec:governing_equations} presents the governing model equations for poroelastic media with deforming and propagating fractures. The discretisation schemes and numerical solution strategy is presented in Section \ref{sec:discretisation}. Section \ref{sec:results} presents simulation results, before concluding remarks are given in Section \ref{sec:conclusion}.

\section{Governing equations}\label{sec:governing_equations}
The conceptual model is based on explicit and conforming representation of fractures in the porous medium. Two modes of fracture deformation are considered: Deformation with fixed transverse extension governed by contact mechanics relations and deformation through irreversible fracture propagation. 
We also impose conservation of mass and energy in matrix and fractures and momentum balance in the matrix.

\subsection{Geometrical representation of fractured porous media}
The model and governing equations are posed in a mixed-dimensional framework arising from considering fractures as lower-dimensional objects. Hence, in a three-dimensional domain, fractures are represented as two-dimensional surfaces, and in a two-dimansional domain, they are one-dimensional lines. In a \nd-dimensional domain, we denote the matrix subdomain by \domain[h] and fractures are represented by subdomains \domain[l] of dimension $\nd-1$. 
The matrix and fractures are connected by interfaces denoted by \interface[j], with the subscript pair $j, k$ used to indicate the two interfaces on either side of a fracture, see Fig.~\ref{fig:tip_coordinate_system}. The boundary of \domain[i] is denoted by \boundary{i}{}, and the internal part of it corresponding to \interface[j] is \boundary{i}{j}.

We also use subscripts $i$, $h$ and $l$ to identify the domain of the primary variables, which are displacement, pressure, temperature, contact traction (\displacement, \pressure, \temperature and \traction). 
Similarly, subscript $j$ denotes the four interface variables defined in Sections \ref{sec:contact_mechanics_equations} and \ref{sec:interface_equations}.
The subscripts are suppressed when context allows, as are the subscripts $f$ and $s$ denoting fluid and solid, respectively.

To model fracture deformation, it is necessary to decompose a vector into its  normal and tangential components relative to a fracture.
The fracture normal is defined to equal the outwards normal \normal[h] on the $j$ side, i.e.\   $\normal[l] =  \normal[h]|_{\boundary{h}{j}}$. 
A vector $\vectorFont{\iota}_l$ may now be decomposed as
\begin{equation}\label{eq:fracture_decomposition}
    \iota_{n} = \vectorFont{\iota}_l \cdot \normal[l] \text{ and }   \vectorFont{\iota}_{\tau} = \vectorFont{\iota}_l - i_{n} \normal[l],
\end{equation}

\noindent
where subscripts $n$ and $\tau$ denote the normal and tangential direction, respectively.



\subsection{Contact mechanics for fracture slip and opening}\label{sec:contact_mechanics_equations}
\begin{figure}
  \includegraphics[width=.49\textwidth]{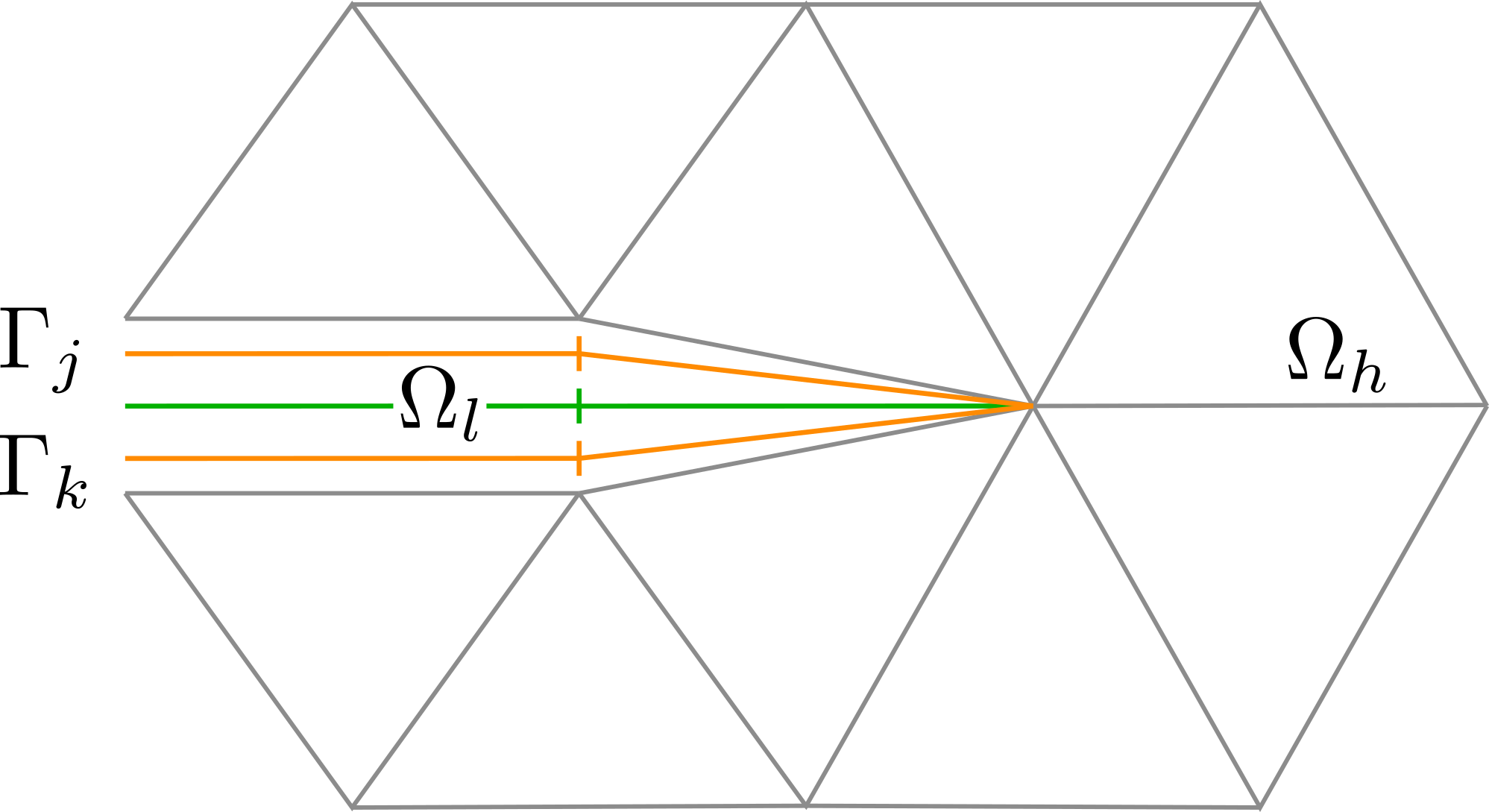}
  \hfill
  \includegraphics[width=.44\textwidth]{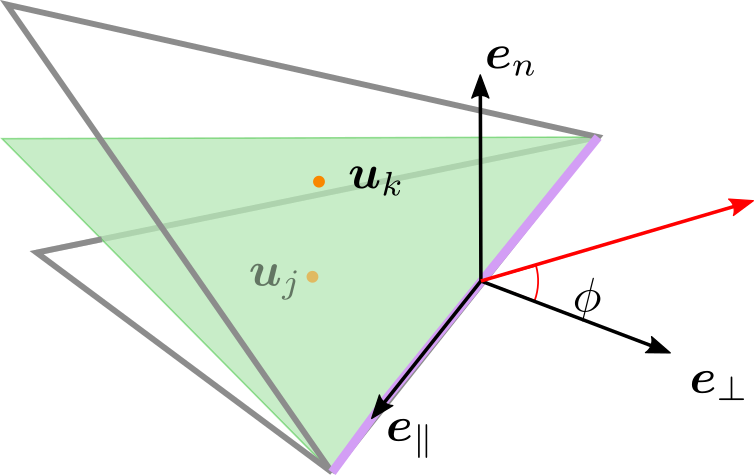}
\caption{Left: A two-dimensional matrix domain \domain[h] and a one-dimensional fracture \domain[l] connected by interfaces, all gridded in a conforming way. Right: Local coordinate system at tip of a two-dimensional fracture. The face of the tip cell is shown in purple and the interface cell centres on the $j$ and $k$ sides are shown as orange dots. The red propagation vector forms an angle \propagationAngle with \ePerp.
The separation between fracture, interfaces and (to the left) matrix faces is for visualisation purposes only, in the model, all coincide geometrically.}
\label{fig:tip_coordinate_system}     
\end{figure}

The contact mechanics relations are a traction balance between the two fracture surfaces and a nonpenetration condition, complemented by a Coulomb friction law governing the relative displacement when the surfaces are in contact. These relations are formulated in the displacement jump \jump{\displacement[]} and the contact traction \traction[l].
The higher-dimensional THM traction, $\stress[h] \cdot \normal[h]$, is balanced by the contact traction and the fracture pressure on the two interfaces:
  \begin{gather}\label{eq:interface_traction_balance}
  \begin{aligned}
  (\traction[l]-\pressure[l] \identity\cdot \normal[l])|_{\Omega_l \cup \Gamma_j} &= \stress[h] \cdot \normal[h]|_{\partial \Omega_h \cup \Gamma_j} 
  , \\
  (\traction[l]-\pressure[l] \identity\cdot \normal[l])|_{\Omega_l \cup \Gamma_k} &= -  \stress[h] \cdot \normal[h]|_{\partial \Omega_h \cup \Gamma_k}.
    \end{aligned}
  \end{gather}
Here the notation indicating that the variable is taken at the interface $Gamma_j$ or $\Gamma_k$ should be interpreted as the extension and projection of this variable to the respecitve interface.
The displacement jump  over the fracture is defined as 
\begin{gather}\label{eq:displacement_jump_definition}
  \begin{aligned}
\jump{ \displacement[l]} = \displacement[k] - \displacement[j],  
  \end{aligned}
\end{gather}
with \displacement[j] and \displacement[k] denoting displacement at \interface[j] and \interface[k], cf.\ Fig.\ \ref{fig:tip_coordinate_system}.
The gap function \gap is defined as the normal distance between the fracture surfaces when these are in mechanical contact.
Following Stefansson \cite{stefansson2020thm}, we set 
\begin{equation}\label{eq:gap_of_displacement}
  \begin{aligned}
  \gap= \tan(\dilationAngle)\norm{\jump{\displacement}_{\tangential}},
  \end{aligned}
  \end{equation}
with \dilationAngle denoting the dilation angle \cite{barton1976shear}, thus accounting for shear dilation of the fracture resulting from tangential displacement $\jump{\displacement}_{\tangential}$ of the rough fracture surfaces.

Given that fracture surface interpenetration and positive normal contact traction are prohibited, the following conditions have to be fulfilled:
\begin{align}\label{eq:fracture_nonpenetration}
\begin{array}{ r l }
                \jump{\displacement}_n - \gap & \geq 0,  \\
                \normalTraction (\jump{\displacement}_n-\gap) & = 0,  \\
                \normalTraction & \leq 0.
  \end{array}  
\end{align}
Hence, when a fracture is mechanically open and there is no mechanical contact across the fracture, the normal contact force, $\lambda_n$, is zero. 

The friction law is imposed by enforcing
\begin{align}\label{eq:fracture_Coulomb}
     \begin{array}{ r l }
                \norm{\traction[\tau]}& \leq -F\normalTraction , \\
                 \norm{\traction[\tau]}& < -\frictionCoefficient \normalTraction \rightarrow \incrementTangentialDisplacement = 0 ,\\
                 \norm{\traction[\tau]}& = -\frictionCoefficient \normalTraction \rightarrow \exists \,\zeta \in \mathbb{R^+}:   \incrementTangentialDisplacement = \zeta\traction[\tau], 
  \end{array}  
\end{align}
where \frictionCoefficient and \incrementTangentialDisplacement denote the friction coefficient and the tangential (shear) displacement increment, respectively. For simplicity, we consider a constant coefficient of friction in this work.

\subsection{Fracture propagation}\label{sec:propagation_equations}
Fracture propagation occurs when the potential energy released by the extension exceeds the energy required to separate the fracture surfaces by breaking atomic bonds \cite{griffith1921rupture}.
Using concepts of linear elastic fracture mechanics, we evaluate propagation  based on computation of stress intensity factors (SIFs). The SIFs are computed directly from the displacement jump in the vicinity of the fracture tip using the variant of the displacement correlation method \cite{nejati2015}. 

Referring to Fig.~\ref{fig:tip_coordinate_system}, the local geometry at a fracture tip is described using a coordinate system given by the orthogonal basis vectors \ePerp, \eN and \ePar (or \ePerp and \eN if $\nd=2$). 
We set $\eN=\normal[l]$, and the tangential (\tangential) vectors \ePerp and \ePar are respectively perpendicular and parallel to \boundary{l}{} at \boundary{l}{} (see Fig.~\ref{fig:tip_coordinate_system}).
By use of the components of the displacement jump in the local coordinate system, the displacement correlation method gives the three SIFs 
\begin{gather}\label{eq:SIFs_DC}
\begin{aligned}
\SIF{I} &= \sqrt{\frac{2\pi}{\DCr}}\bigg(\frac{\shearModulus}{\kolosov+1}\jump{\displacement}_n\bigg), \\
\SIF{II} &= \sqrt{\frac{2\pi}{\DCr}}\bigg(\frac{\shearModulus}{\kolosov+1}\jump{\displacement}_{\perp}\bigg) ,\\
\SIF{III} &= \sqrt{\frac{2\pi}{\DCr}}\bigg(\frac{\shearModulus}{4}\jump{\displacement}_{\parallel}\bigg).
\end{aligned} 
\end{gather}
Here, \shearModulus denotes the shear modulus and $\kolosov=3-4\poisson$ the Kolosov constant, with \poisson being the Poisson ratio. \DCr is the distance from the fracture tip to the point at which the displacement jump is evaluated.
The three stress intensity factors are related to tensile ($\SIF{I}$), shear ($\SIF{II}$) and torsional ($\SIF{III}$) forces.

As stated in the introduction, we limit ourselves to tensile fracture in this work, and ignore contributions from ($\SIF{II}$) and ($\SIF{III}$). 
With this assumption,  a tip propagates if the computed mode I factor exceeds a critical value, \begin{gather}\label{eq:propagation_criterion}
\begin{aligned}
\SIF{I} \geq \SIF{Ic},
\end{aligned} 
\end{gather}
and the propagation angle \propagationAngle illustrated in Fig.~\ref{fig:tip_coordinate_system} is zero. 
Criteria for more sophisticated mixed-mode propagation, which can be highly relevant for subsurface applications, are reviewed by Richard et al. \cite{richard2005fracturepropagation}.


\subsection{Fracture mass and energy balance}\label{sec:fracture_em}
The thickness of a dimensionally reduced fracture is represented by the aperture, which changes as the domain deforms according to
 \begin{equation}
     \aperture = \aperture[res] + \jump{\displacement}_n,
 \end{equation}
with \aperture[res] denoting the residual hydraulic aperture in the undeformed state representing the effect of small-scale roughness of the two fracture surfaces. The tangential fracture permeability $\permeability[l]$ is chosen to depend on aperture by the nonlinear relationship $\permeability[l]=a^2/12$, which corresponds to setting the hydraulic aperture of the fracture equal to \aperture \cite{zimmerman1996cubic}.

On the assumption that the fractures are completely filled with fluid, the parameters of this subsection equal that of the fluid.
We assume single-phase flow according to Darcy's law:
\begin{gather}\label{eq:darcy}
\begin{aligned}
 \fluidFlux = -\frac{\permeability}{\viscosity} \left(\nabla \pressure - \density[]\gravityVector\right),
\end{aligned}
\end{gather} 
where \viscosity, \density and \gravityVector denote viscosity, density and the gravity acceleration.
The total heat flux may be split into continuum scale heat diffusion modelled by Fourier's law and advection along the fluid flow field:
\begin{gather}\label{eq:heat_fluxes}
\begin{aligned}
 \conductiveHeatFlux &= -\heatConductivity \nabla \temperature, \\
 \advectiveHeatFlux[] &=  \density \heatCapacity \temperature \fluidFlux,
\end{aligned}
\end{gather}
where \heatConductivity and  \heatCapacity denote thermal conductivity and heat capacity, respectively.

Following Stefansson et al. \cite{stefansson2020thm} (see also Brun et al. \cite{brun2018upscaling} and Coussy \cite{coussy2004poromechanics}), balance of mass for a fracture \domain[l] reads 
\begin{gather}\label{eq:reduced_mass_balance_darcy}
\begin{aligned}
\aperture\left(\compressibility\dd{\pressure}{\timet} - \thermalExpansion\dd{\temperature}{\timet}\right) +\dd{\aperture}{\timet}
-  \divergence\left(\aperture\frac{\permeability}{\viscosity}\left(\nabla\pressure-\density\gravityVector\right)\right) - \sum_{j \in \higherSet[l]}\interfaceFluidFlux[j] =  \aperture \source[\pressure],
\end{aligned}
\end{gather}
with \compressibility, \thermalExpansion and \source[\pressure] denoting compressibility, thermal expansion coefficient and  a fluid source or sink term.

Next, assuming local thermal equilibrium between fluid and solid, neglecting viscous dissipation and linearising \cite{stefansson2020thm}, the energy balance is
\begin{gather}\label{eq:entropy_balance_linearised_simplified}
\begin{aligned}
 &\frac{\heatCapacity\density}{\temperature[0]}(\temperature-\temperature[0])\dd{\aperture}{\timet} 
+ \frac{\heatCapacity\density}{\temperature[0]}
\aperture\dd{\temperature}{\timet} 
- \thermalExpansion \aperture\dd{\pressure}{\timet} 
+ \divergence\left[ \aperture\left( \frac{\heatCapacity\density}{\temperature[0]}(\temperature-\temperature[0]) \fluidFlux- \frac{\heatConductivity}{\temperature[0]} \nabla \temperature \right)\right] \\
&- \sum_{j \in \higherSet[l]}\frac{\interfaceConductiveHeatFlux[j]}{\temperature[0]} +\frac{\interfaceAdvectiveHeatFlux[j]}{\temperature[0]} 
=  \aperture\source[\temperature],
\end{aligned}
\end{gather}
where we assume thermal sources and sinks to satisfy $\aperture\source[\temperature]=\aperture \source[\pressure] \frac{\heatCapacity\density}{\temperature[0]}(\temperature-\temperature[0])$ and \temperature[0] denotes a reference temperature. 
In Eqs.~\eqref{eq:reduced_mass_balance_darcy} and \eqref{eq:entropy_balance_linearised_simplified} the last terms on the right hand sides represents the fluxes from matrix to fractures, which are defined in Section \ref{sec:interface_equations}.

In deriving these equations, the following equations of state are assumed \cite{coussy2004poromechanics} for density
\begin{gather}\label{eq:density_eos}
\begin{aligned}
\density &= \density[0] \exp[\compressibility (\pressure-\pressure[0]) - \thermalExpansion (\temperature-\temperature[0])]
\end{aligned}
\end{gather}
and entropy
\begin{gather}\label{eq:entropy_eos}
\begin{aligned}
\entropy - \entropy[0] &= -\thermalExpansion \frac{\pressure-\pressure[0]}{\density} + \frac{\heatCapacity}{\temperature[0]}(\temperature - \temperature[0]).
\end{aligned}
\end{gather}

\subsection{Matrix thermo-poroelasticity, energy and mass balance}\label{sec:thm_equations}
The following section presents the balance equations and constitutive relations for the matrix problem. 
The model resembles that of the previous section, with the addition of a momentum balance equation for the thermo-poroelastic medium, yielding three balance equations for \domain[h]. For details on the derivations of the equations, we  again refer to Coussy \cite{coussy2004poromechanics} and Brun et al. \cite{brun2018upscaling}.
We first define the following effective parameters \cite{cheng2016poroelasticity}, arising through the assumption of local thermal equilibrium:
\begin{gather}\label{eq:effective}
\begin{aligned}
 \heatConductivity[\eff]        &= \porosity\heatConductivity[f] + (1-\porosity)\heatConductivity[s],\\
  (\density \heatCapacity)_\eff &= \porosity\density[f]\heatCapacity[f] + (1-\porosity)\density[s]\heatCapacity[s],\\
   \thermalExpansion[\eff]      &= \porosity\thermalExpansion[f] + (\biotAlpha-\porosity)\thermalExpansion[s].
  \end{aligned}
\end{gather}
\porosity and  \biotAlpha denote porosity and the Biot coefficient, respectively.

Neglecting inertial terms, the momentum balance is
\begin{gather}\label{eq:momentum_balance}
\begin{aligned}
\divergence \stress = \vectorSource[\displacement], 
\end{aligned} 
\end{gather}
with \vectorSource[\displacement] denoting body forces and the linearly thermo-poroelastic stress tensor related to the primary variables  by an extended Hooke's law 
\begin{gather}\label{eq:mechanics_stress}
  \begin{aligned}
    \stress - \stress[0] = \frac{\stiffness}{2}(\nabla\displacement + \nabla\displacement^T)- \biotAlpha (\pressure-\pressure[0]) \identity - \thermalExpansion[\eff]  \bulkModulus (\temperature - \temperature[0]) \identity.
  \end{aligned}
\end{gather}
The mass balance equation reads
\begin{gather}\label{eq:mass_balance}
\begin{aligned}
\left(\porosity \compressibility+\frac{\biotAlpha-\porosity}{\bulkModulus}\right) \dd{\pressure}{\timet} + \biotAlpha \dd{(\divergence\displacement)}{\timet}-\thermalExpansion[f]\dd{\temperature}{\timet} + \divergence \left(\frac{\permeability}{\viscosity}\left(\nabla\pressure-\density\gravityVector\right)\right) = \source[\pressure],
\end{aligned}
\end{gather}
while the energy balance is
\begin{gather}\label{eq:heat_balance}
\begin{aligned}
\frac{(\density \heatCapacity)_\eff}{\temperature[0]} \dd{\temperature}{\timet} + \thermalExpansion[s] \bulkModulus  \dd{(\divergence\displacement)}{\timet}-\thermalExpansion[f]  \dd{\pressure}{\timet} + \divergence\left( \frac{\heatCapacity\density}{\temperature[0]}(\temperature-\temperature[0]) \fluidFlux- \frac{\heatConductivity}{\temperature[0]} \nabla \temperature \right)  = \source[\temperature].
\end{aligned}
\end{gather}

On $\Omega_h \cup \Gamma_j$, the following internal boundary conditions ensure coupling from \domain[h] to the interface variables on \interface[j]:
\begin{equation}\label{eq:internal_boundary_conditions}
  \begin{aligned}
    \displacement[h] &=\displacement[j], \\
	\fluidFlux[h] \cdot \normal[h]&=\interfaceFluidFlux,  \\
    \conductiveHeatFlux[h] \cdot \normal[h]&=\interfaceConductiveHeatFlux[j],\\
    \advectiveHeatFlux[h] \cdot \normal[h]&=\interfaceAdvectiveHeatFlux[j].
	\end{aligned}
  \end{equation}
The conservation equations are complemented by appropriate boundary conditions on the domain bondary. This applies to both the matrix and fracture domains.

\subsection{Interface fluxes between fractures and matrix}\label{sec:interface_equations}
Interface flux relations close the mixed-dimensional system of mass and energy balance equations \cite{Martin2005,jaffre2011interdimupwind}:
\begin{gather}\label{eq:interface_fluxes}
    \begin{aligned}
    \interfaceFluidFlux &= - \frac{\permeability[j]}{\viscosity}\left( \frac{2}{\aperture[l]} \left(  \pressure[l]|_{\Omega_l \cup \Gamma_j} -  \pressure[h]|_{\partial \Omega_h \cup \Gamma_j}  \right) -\density[l] \gravityVector \cdot \normal[h] \right),  \\
    \interfaceConductiveHeatFlux &= - \heatConductivity[j]\frac{2}{\aperture[l]} ( \temperature[l]|_{\Omega_l \cup \Gamma_j} - \temperature[h]|_{\partial \Omega_h \cup \Gamma_j}),  \\
    \interfaceAdvectiveHeatFlux &= \left\{ \begin{array}{ l l }
      \interfaceFluidFlux \density[h]\heatCapacity[h]\temperature[h] & \text{ if } \interfaceFluidFlux>0  \\
      \interfaceFluidFlux  \density[l]\heatCapacity[l]\temperature[l] & \text{ if } \interfaceFluidFlux\leq0 
  \end{array} \right. .
     \end{aligned}
\end{gather} 
We set the normal permeability and thermal conductivity equal to their tangential counterparts, i.e.\ $\permeability[j]=\permeability[l]$ and $\heatConductivity[j]=\heatConductivity[l]$.

\section{Discretisation and solution strategy}\label{sec:discretisation}
Discretisation of the governing equations entails devising discrete representation of the conservation equations and of the contact mechanics relations on existing fractures.
Moreover, when the propagation criteria are met, the fracture geometry must be modified and the discretisations updated accordingly.

We make the following assumptions on the computational grid:
Grids for the subdomains \domain[h] and \domain[l] and the interface \interface[j] are constructed so that faces on $\partial_j\Omega_i$ match with cells in \interface[j] and \domain[l].
We make no assumptions on the cell types; for the simulations presented in Section \ref{sec:results} we mainly use Cartesian grids as these are most easily fit to a known, straight propagation path, but also consider simplex cells for one simulation.

\subsection{Spatial discretisation}
Pressure and temperature are represented by their cell centre values in \domain[h] and \domain[l], as is the displacement in \domain[h] and contact force in \domain[l].
The discrete primary variables on \interface[] are displacements, mass flux and advective and diffusive heat fluxes. 

\subsubsection{Contact mechanics}
The non-linear contact mechanics problem is represented by an active set approach implemented as a semi-smooth Newton method following \cite{hueber2008paper,berge2020hmdfm}.
The treatment of Eqs. \eqref{eq:fracture_nonpenetration} and \eqref{eq:fracture_Coulomb} depends on whether the previous iterates were in an open, sticking or gliding state, with the states evaluated cell-wise in \domain[l]. 
Equation \eqref{eq:interface_traction_balance} is discretised by relating the cell centre pressures and contact force in \domain[l] to the discrete traction on \interface[].

\subsubsection{Discretisation of balance equations}
For \domain[h], the stress term in \eqref{eq:momentum_balance} and the diffusive fluxes in \eqref{eq:mass_balance} and \eqref{eq:heat_balance} are all discretised with a family of finite volume multi-point approximations termed MPxA \cite{Aavatsmark2002,nordbotten2016biot,nordbotten2020mpxa}. 
The methods construct discrete representations of the constitutive relations, Hook's, Darcy's and Fourier's law, in terms of the cell centre variables.
These relations are used to enforce conservation of THM traction, mass and (diffusive) heat flux over the cell faces.
For faces on the fracture surfaces, the discrete traction enters the contact mechanics discretisation described above.
The full heat flux is given by the sum of the discrete Fourier's law and the advective flux, where the latter is discretised by a single-point upstream method.
For further information on the MPxA methods, we refer to \cite{nordbotten2020mpxa}. 

In \domain[l], Eqs. \eqref{eq:reduced_mass_balance_darcy} and \eqref{eq:entropy_balance_linearised_simplified} are discretised analogously to the corresponding terms in \domain[h].
Finally, fluxes over \interface are computed from discrete versions of Eqs. \eqref{eq:interface_fluxes}.



\subsubsection{Solution of non-linear system}
The discretised system of equations is solved by Newton's method, with the terms from the contact conditions handled by a semi-smooth approach following \cite{hueber2008paper,berge2020hmdfm}.
The termination criterion for the Newton iterations considers the residuals and updates of each of the primary variables \displacement[h], \displacement[j], \pressure and \temperature.
Within each non-linear iteration, the linearised system is solved using a direct sparse solver \cite{umfpack}. While simple, this approach is memory intensive and puts practical constraints on mesh resolution, in particular for three-dimensional problems.
A more scalable method would involve iterative solvers with block preconditioners for the THM components of the linear system \cite{both2019gradient}, with a tailored treatment of the contact conditions \cite{franceschini2019block}.

\subsection{Solution algorithm}
The temporal derivatives are discretised by a backward Euler scheme, and the THM contact mechanics problem is solved monolithically, using implicit in time evaluation of all spatial derivatives. When the non-linear solver has converged, we proceed to fracture propagation evaluation.

\begin{figure}
  \includegraphics[width=\textwidth]{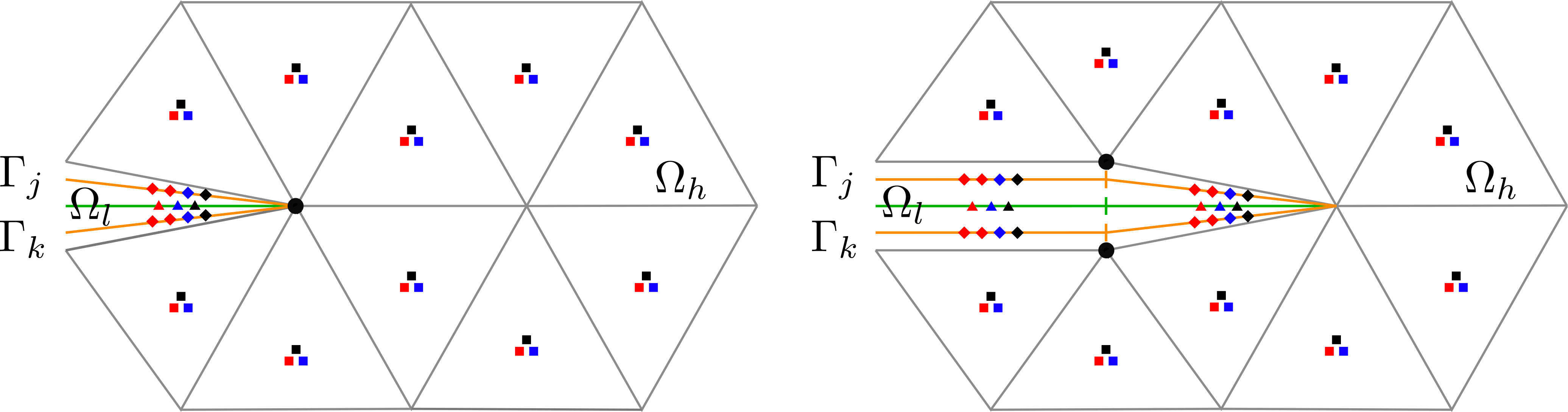}
\caption{Example grid and unknowns before (left) and after (right) propagation. Unknown shapes reflect the subdomain or interface where they are defined, whereas the colours red, blue and black correspond to heat, mass and deformation, respectively. One vertex (black circle) and the face along which propagation occurs are duplicated as part of the geometry update. The different types of domains are separated for illustration purposes. Figure adapted from \cite{stefansson2020thm}.}
\label{fig:propagation_grid} 
\end{figure}

Stress intensity factors and the fracture propagation criterion are evaluated for each fracture tip faces using Eqs.~\eqref{eq:SIFs_DC} and~\eqref{eq:propagation_criterion}. 
The displacement jump is evaluated at the neighbouring cell of the tip face, i.e. \DCr is the distance between the centre of the tip face and the cell. The fracture is restricted to grow along faces of the matrix grid, no computation of the propagation length is performed. 
The geometrical update for each identified face now entails i) duplicating the face for the matrix grid, ii) adding a cell in the fracture grid and iii) adding one cell for each of the two interfaces. The three new cells all coincide geometrically with the chosen face, see Fig.~\ref{fig:propagation_grid}.
Once new cells and faces have been added, connectivity information is updated both within subdomains and between the subdomains and the interface.

Variables are initialised in the new cells using the reference values \pressure[0] and \temperature[0], and new apertures are set to \aperture[res].
This in effect adds mass to the system, cf. Eq.~\eqref{eq:density_eos}.
To compensate, we prescribe an additional term on the right hand side of equation \eqref{eq:reduced_mass_balance_darcy} equal to $-\aperture[res]/dt$ in newly formed fracture cells the subsequent time step, with $dt$ denoting time step size. 
Since Eq.~\eqref{eq:entropy_balance_linearised_simplified} is derived by considering $\entropy-\entropy[0]$, Eq.~\eqref{eq:entropy_eos} implies that no right-hand side term arises with the chosen initialisation values.

Before the simulation proceeds to the next time step, all terms are rediscretised to account for modifications of the grids. 
This can be done locally, i.e. only for the faces and cells where the discretisation is affected by the grid update.


\section{Simulation results}\label{sec:results}
The results presented in this section serve to first verify the computational approach, and then to show application to two subsurface cases involving THM processes and fracture propagation. 
The PorePy toolbox \cite{keilegavlen2020porepy,keilegavlen2020porepyv130} was used for all simulations and run scripts for geometry and parameter setup etc. are available on GitHub \cite{stefansson2020sourcecode_propagation}. All parameters not specified in the text are listed in Table \ref{tab:parameters}.

\subsection{Verification}
The verification of the computational approach entails first a test of the numerical stress intensity factors and next a convergence test of the fracture propagation speed.

\subsubsection{Example 1: Stress intensity factors}\label{sec:ex1}
To verify the SIF computation, we consider an analytical solution, first derived by Sneddon \cite{sneddon1946analytical}, for a single crack in an infinite medium with uniform internal pressure on the fracture surfaces. 
Boundary conditions for the finite simulation domain are computed using the boundary element method following Keilegavlen et al. \cite{keilegavlen2020porepy}, who also presents a thorough convergence study for the aperture using PorePy.
Herein, we compare the SIFs as computed by the displacement correlation method to the analytical solution
\begin{equation}
    \begin{aligned}
        \SIF{I, an}&= \pressure[f] \sqrt{l \pi},\\ 
        \SIF{II, an}&=0.
    \end{aligned}
\end{equation}
Here, \pressure[f] denotes the internal pressure on the fracture and $l$ denotes fracture length.
\begin{figure}
\centering
\begin{subfigure}{0.5\textwidth}
    \includegraphics[width=\textwidth]{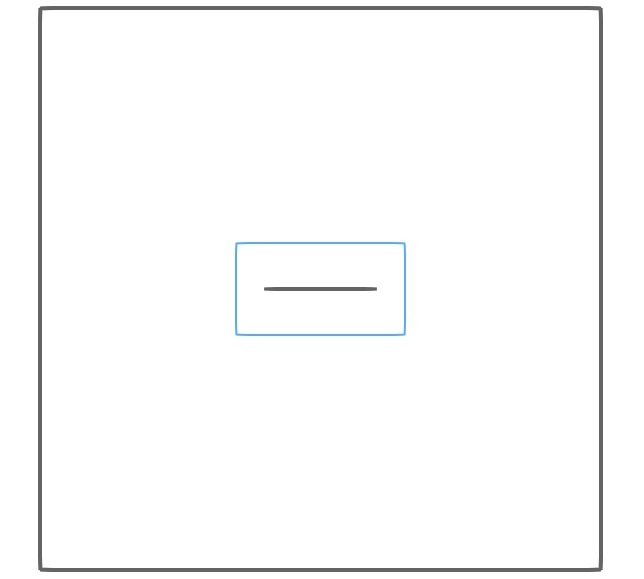}
\end{subfigure}
\hfill
\begin{subfigure}{0.45\textwidth}
\includegraphics[width=\textwidth]{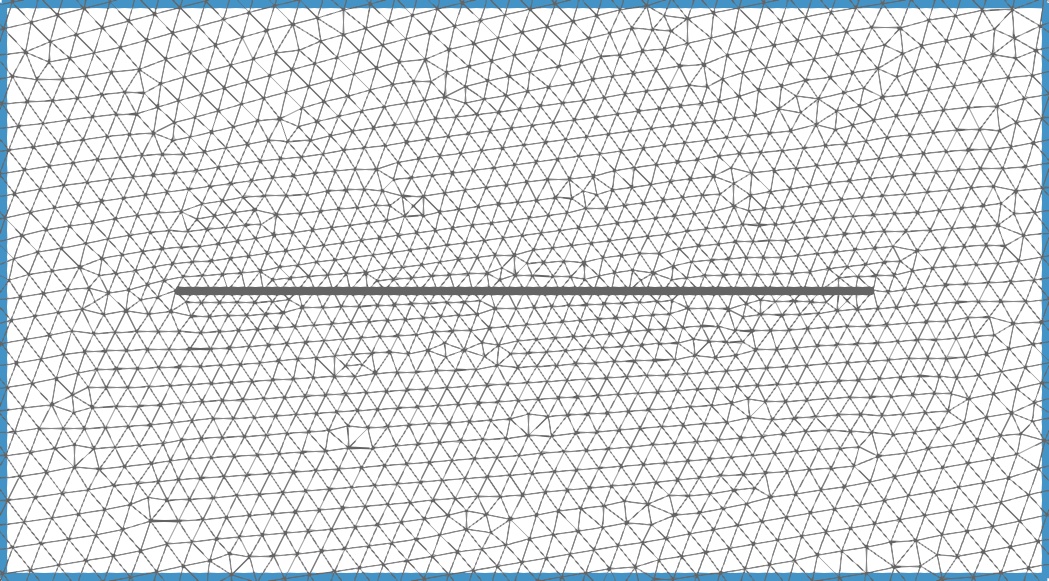} 
\end{subfigure}
\caption{Example 1: Domain geometry (left) and close-up around the fracture showing the matrix mesh for the mesh with $h=0.25$ (right).
The blue box shows the location of the close-up. 
}
\label{fig:ex1_geometry_and_mesh}     
\end{figure}
We use a square domain of side length \SI{50}{\metre}, $l=\SI{10}{\metre}$ and $\pressure[f]=\SI{1e-4}{\pascal}$. We consider a sequence of four grids, the finest of which is shown in Fig.~\ref{fig:ex1_geometry_and_mesh}. To probe the method for different material parameters, we also use four different Poisson ratios.
Based on displacement solutions on each grid, SIFs are estimated and the normalised $L^2$ type errors are computed as
\begin{equation}\label{eq:sif_errors}
    \begin{aligned}
       E_{I} &=  \frac{\left[\sum_{j=1}^2\left(\SIF{I,j}-\SIF{I,an}\right)^2\right]^{1/2}}{2 \SIF{I,an}}, \\  
       E_{II} &=  \frac{\left[\sum_{j=1}^2\left(\SIF{II,j}-\SIF{II,an}\right)^2\right]^{1/2}}{2 \SIF{I,an}}, 
    \end{aligned}
\end{equation}
with the $j$ index running over the two fracture tips.
\begin{figure}
    \includegraphics[width=.49\textwidth]{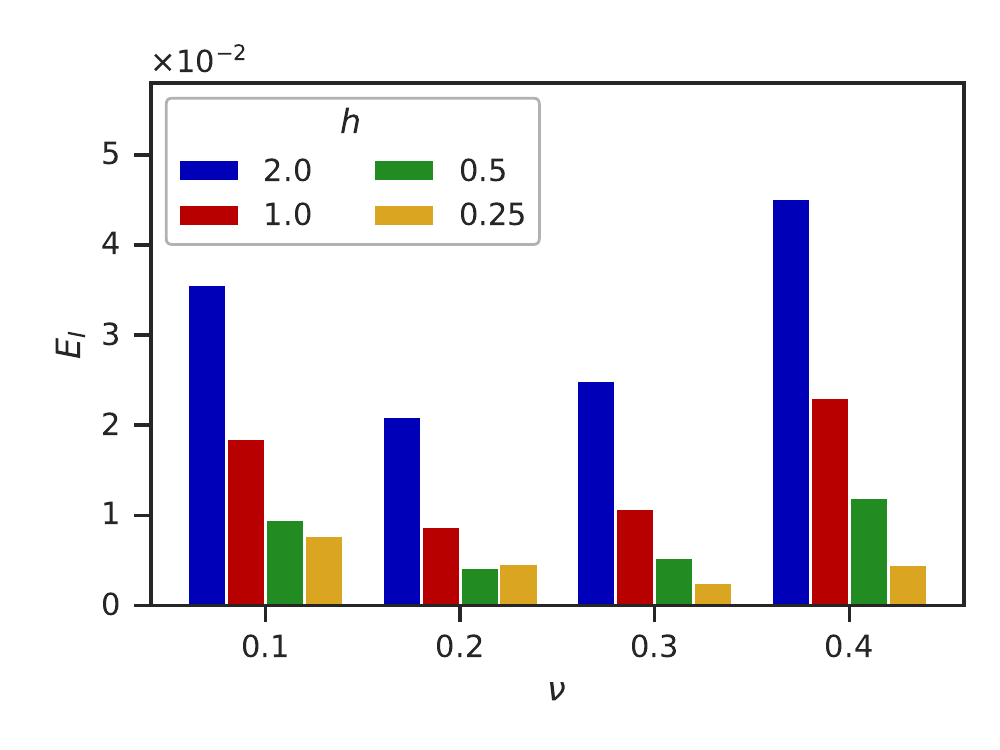}
    \hfill
    \includegraphics[width=.49\textwidth]{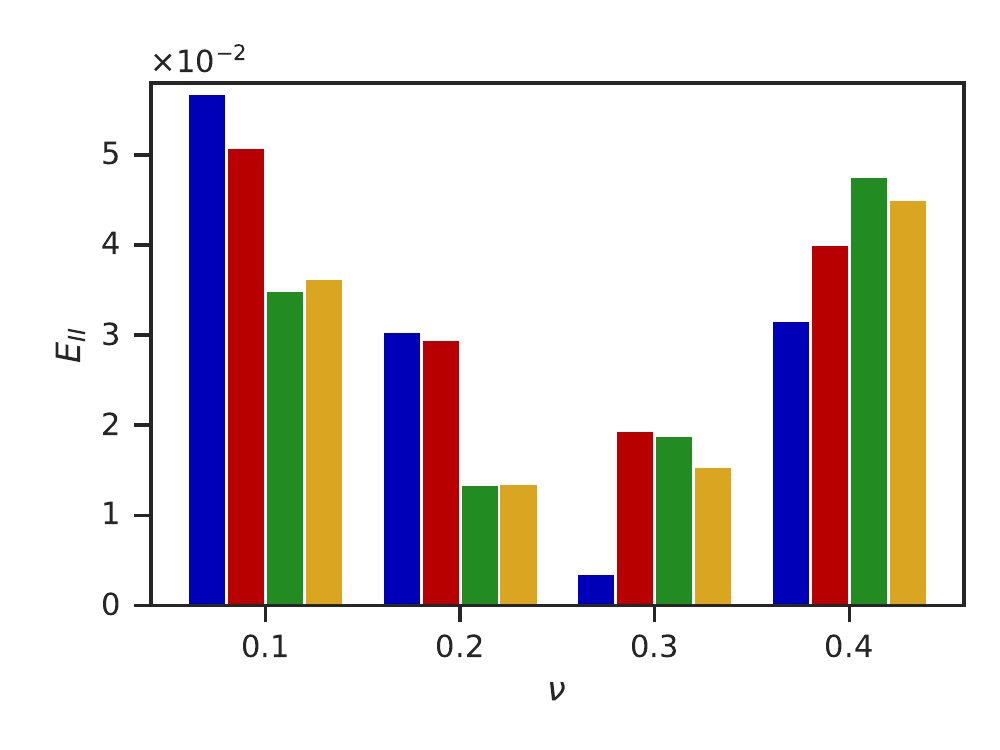} 
\caption{Errors for \SIF{I} (left) and  \SIF{II} (right) computed according to Eq.~\eqref{eq:sif_errors} for different values of the Poisson ratio \poisson and different mesh sizes $h$. }
\label{fig:ex1_errors}     
\end{figure}
The piecewise linear displacement representation of the MPSA discretisation does not capture the stress singularity at the fracture tips. 
Since the SIFs are computed from \jump{\displacement} in these very tip cells, the method does not converge with mesh refinement.
Rather, the results presented in Fig.~\ref{fig:ex1_errors} demonstrate robustness with respect to mesh size and the Poisson ratio \poisson.
While we do not consider \SIF{II} in the subsequent simulations, we also present results demonstrating that the method indeed predicts tensile stresses (i.e. $\SIF{II} \ll \SIF{I}$) for this purely tensile problem.

The results of this test indicate that the MPSA solution can be used to estimate \SIF{I} in tensile problems, and thus form the basis of fracture growth evaluation.

\subsubsection{Example 2: Propagation speed}\label{sec:ex2}
\begin{figure}
\centering
\begin{subfigure}{0.5\textwidth}
    \includegraphics[width=\textwidth]{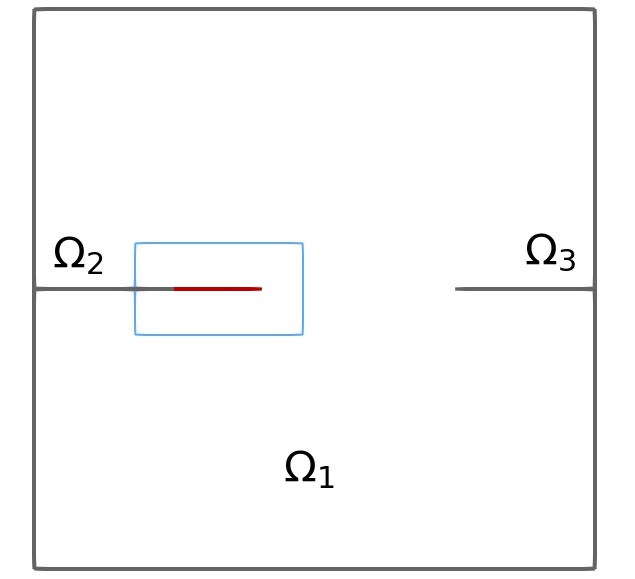}
\end{subfigure}
\hfill
\begin{subfigure}{0.45\textwidth}
    \includegraphics[width=\textwidth]{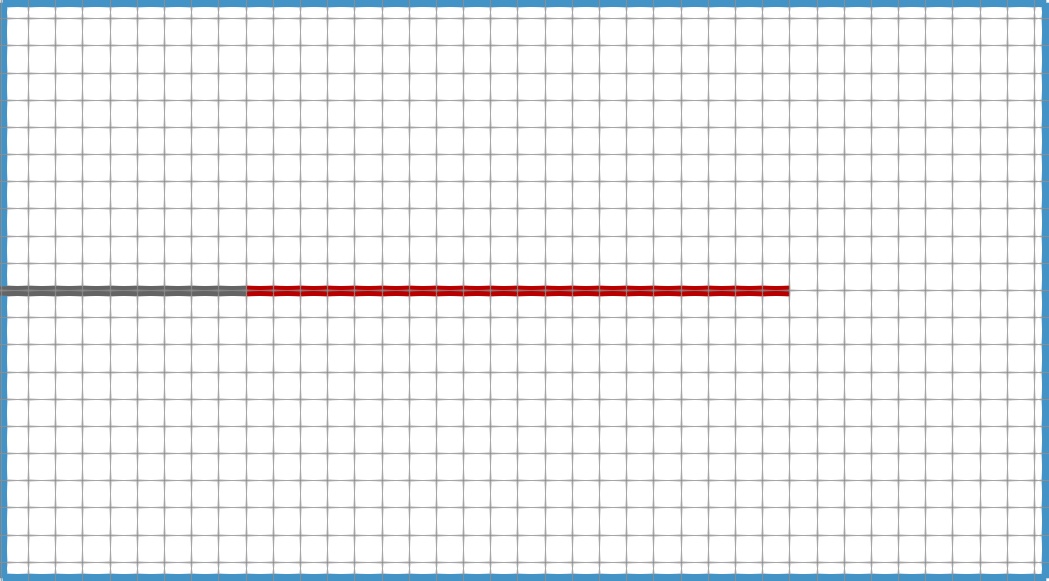} 
\end{subfigure}
\caption{Left: Domain geometry used in \hyperref[sec:ex2]{Example 2}. The grey lines indicate initial geometry, whereas the red line indicates the extension at the end of simulation. The blue box shows the location of the close-up to the right. Right: Close-up around the fracture showing the matrix mesh for refinement $h=1/128$. }
\label{fig:ex2_geometry_and_mesh}     
\end{figure}

We now consider a test case designed to evaluate the simulated propagation speed of a fracture in a tensile regime of stable propagation. 
The unit square domain contains two horizontal fractures \domain[2] and \domain[3] extending 1/4 from the left and right boundary, respectively, see Fig.~\ref{fig:ex2_geometry_and_mesh}. 
The boundary conditions for fluid and heat are no-flow in the matrix and Dirichlet for the fractures, with zero values on the right and $\pressure=\SI{5}{\mega\pascal}$ and $\temperature=\SI{-50}{\kelvin}$ on the left.
Thus, cold fluid flows from left to right, entering the matrix at the right end of \domain[2].
The domain is mechanically fixed at the top and bottom and zero traction is imposed on the left and right boundaries. 
The forces driving propagation are the elevated pressure inside \domain[2] and cooling of the surrounding matrix.


We use four temporal refinement levels and three spatial refinement levels in addition to a highly refined reference solution.
The finest (non-reference) mesh and the final fracture geometry are illustrated in Fig.~\ref{fig:ex2_geometry_and_mesh}.

\begin{figure}
\centering
    \includegraphics[width=.85\textwidth]{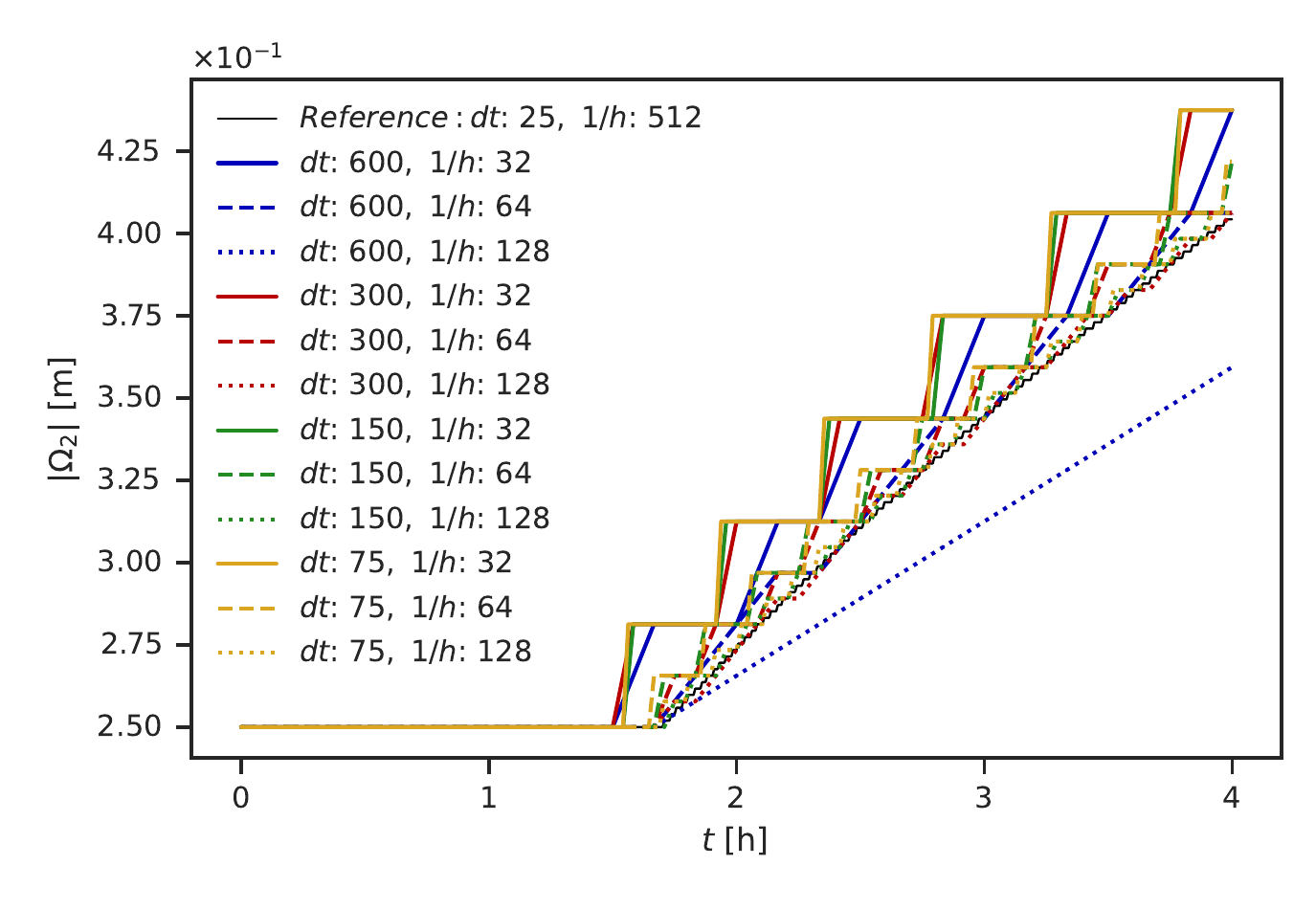}
\caption{\hyperref[sec:ex2]{Example 2}: Size of the propagating fracture \domain[2] vs. time for 13 refinement combinations. Line colours and styles correspond to temporal and spatial discretisation size, respectively.}
\label{fig:ex2_speeds}     
\end{figure}

Figure~\ref{fig:ex2_speeds} shows fracture size plotted against time for all refinement combinations.
With one exception discussed below, the results group according to spatial resolution.
While the propagation speed is fairly constant across all mesh sizes, propagation onset occurs earlier for the coarser meshes.
We attribute this offset to the SIFs being evaluated on the basis of \jump{\displacement}~at the centre of the fracture tip cell. 
The location of this cell centre is closer to the boundary for the coarser meshes, implying shorter travel time for the cooling front. As expected, the plot indicates convergence with mesh refinement.

The outlier is the smallest mesh size combined with the largest time step, for which the propagation speed is notably lower. 
The propagation speed is simply not resolved by the spatio-temporal discretisation, i.e. the propagation speed exceeds $h/dt$. 
In other words: Given a spatial resolution, an upper bound on the time step must be honoured in the explicit type of propagation solution algorithm used herein.

\begin{table}
\caption{Parameters for the simulation examples. For \hyperref[sec:ex1]{Example 1}, only mechanical parameters are relevant.}
\label{tab:parameters}       
\begin{tabular}{lllll}
\hline\noalign{\smallskip}
Parameter & Symbol & Examples & Value & Units   \\
\noalign{\smallskip}\hline\noalign{\smallskip}
Biot coefficient & \biotAlpha & 2-4 & \num{0.8} &\si{-}\\
Friction coefficient & \frictionCoefficient & 2-4 & \num{0.8} &\si{-} \\
Dilation angle & \dilationAngle & 2-4 & \num{3.0} & \si{\degree} \\
Fluid linear thermal expansion  &\thermalExpansion[f] & 2-4 & \num{4e-4}& \si{\per \kelvin} \\
Solid linear thermal expansion  & \thermalExpansion[s] & 2-4 & \num{5e-5}& \si{\per \kelvin} \\
Critical stress intensity factor & \SIF{c} & 2-4 & \num{5e5}& \si{\pascal} \\
Fluid specific heat capacity &\heatCapacity[f] & 2-4 & \num{4.2e3}& \si{\joule\per\kilogram \per\kelvin} \\
Solid specific heat capacity &\heatCapacity[s] & 2-4 & \num{7.9e2}& \si{\joule\per\kilogram \per\kelvin} \\
Fluid thermal conductivity &\heatConductivity[f] & 2-4 & \num{0.6}& \si{\watt\per\meter\per\kelvin} \\
Solid thermal conductivity &\heatConductivity[s] & 2-4 & \num{2.0}& \si{\watt\per\meter\per\kelvin} \\
Reference fluid density &\density[0, f] & 2-4 & \num{1e3}& \si{\kilogram\per\meter\cubed} \\
Reference solid density &\density[0, s] & 2-4 & \num{2.7e3}& \si{\kilogram\per\metre\cubed} \\
Compressibility & \compressibility & 2-4 & \num{4e-10}& \si{\per\pascal} \\
Bulk modulus & \bulkModulus & 1-4 & \num{2.2e10}& \si{\pascal} \\
Poisson ratio & \poisson & 1-4  & \num{0.2} & \si{-}\\ 
Matrix porosity & \porosity & 2-4 & \num{0.05} & \si{-}\\
Matrix permeability & \permeability & 2-3 & \num{1e-14}& \si{\metre\square} \\
Matrix permeability & \permeability & 4 & \num{1e-16}& \si{\metre\square} \\
Viscosity & \viscosity & 2-4 & \num{1e-3}& \si{\pascal\second} \\
Residual aperture & \aperture[res] & 2 & \num{1e-3}& \si{\metre} \\
Residual aperture & \aperture[res] & 3 & \num{3e-4}& \si{\metre} \\
Residual aperture & \aperture[res] & 4 & \num{2.e-3}& \si{\metre} \\
\noalign{\smallskip}\hline
\end{tabular}
\end{table}

\subsection{Applications}
We present two simulations that involve THM processes coupled with fracture propagation.
The first case resembles geothermal energy production, with convection forced by fluid injection and production.
The second case involves natural convection that takes place mainly inside fractures.
In both cases, convection acts to alter thermal stresses and thereby cause fracture propagation.

\subsubsection{Example 3: Thermal fracturing and forced convection}\label{sec:ex3}
\begin{figure}
\centering
\begin{subfigure}{0.5\textwidth}
\includegraphics[width=.99\textwidth]{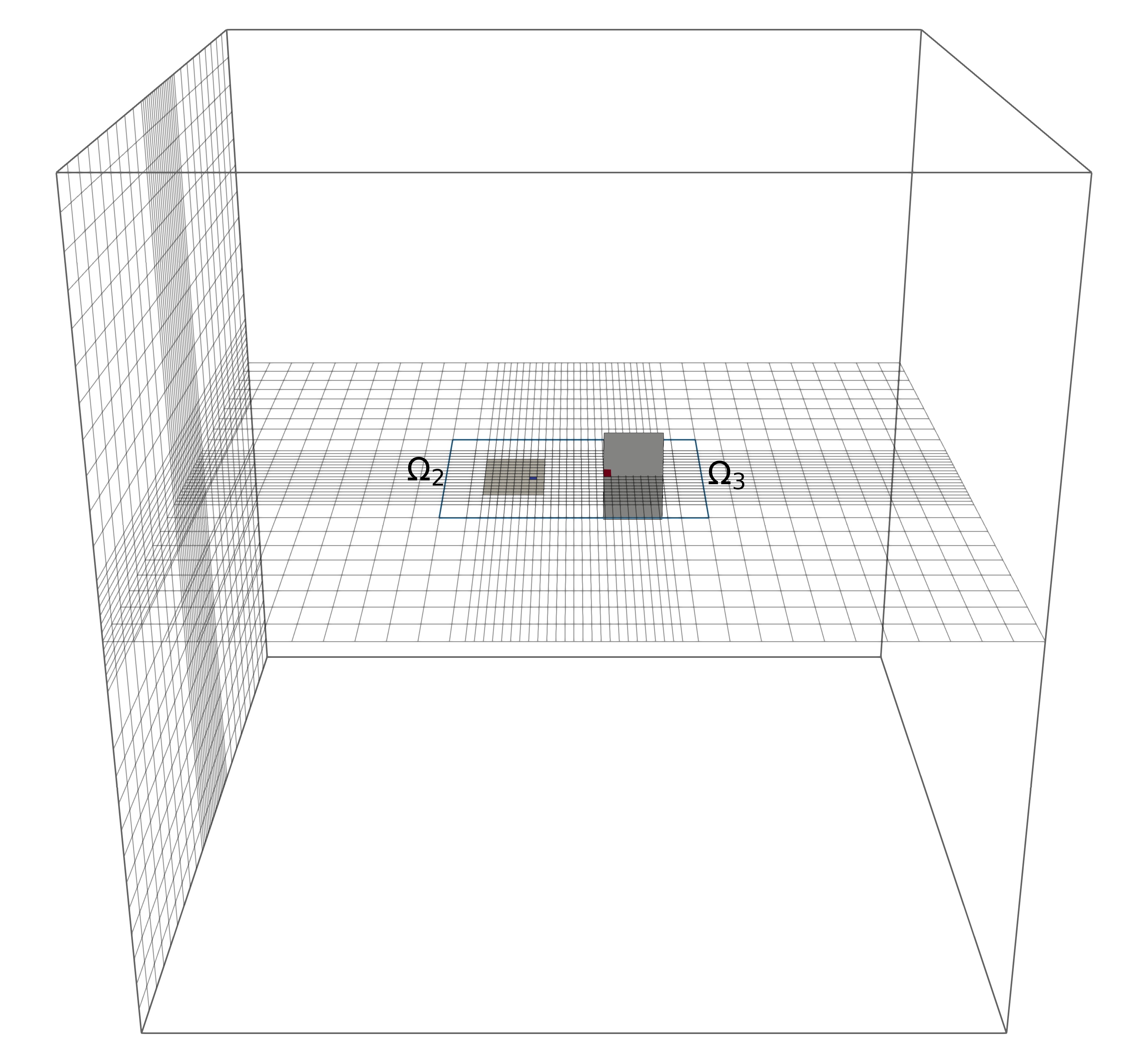}
\end{subfigure}
\hfill
\begin{subfigure}{0.46\textwidth}
        \includegraphics[width=\textwidth]{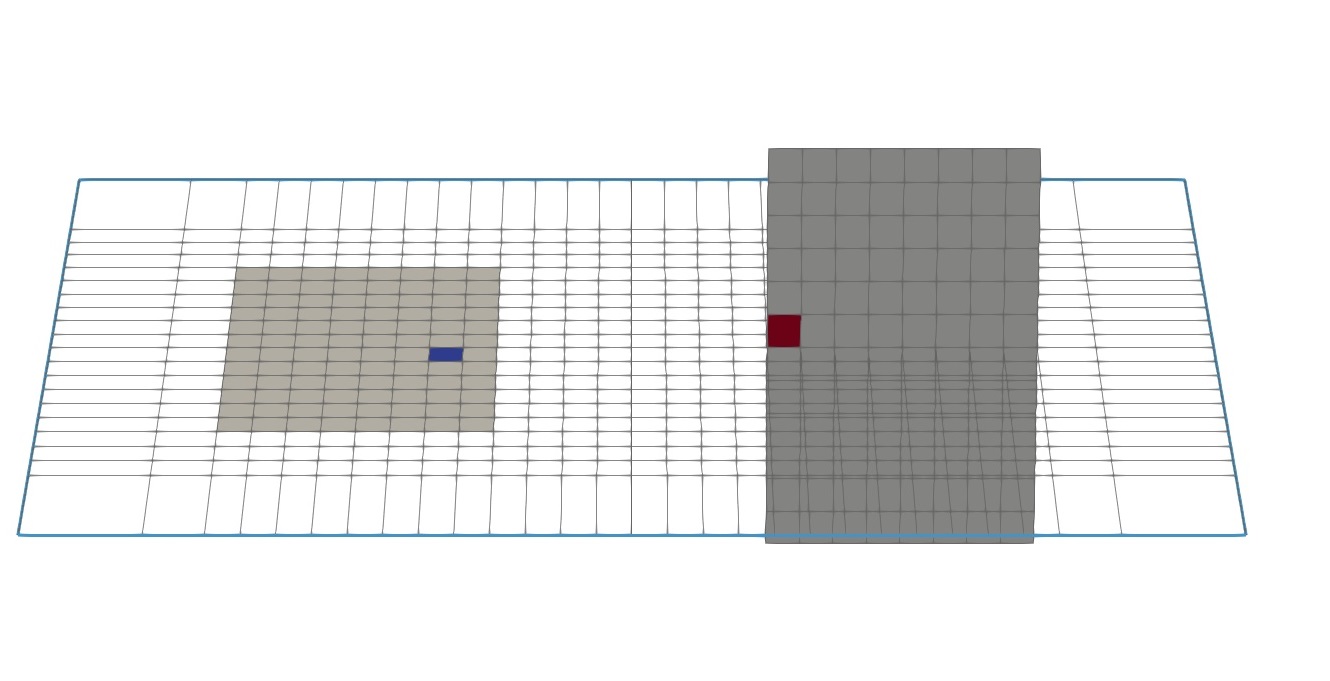}
\end{subfigure}    
\caption{\hyperref[sec:ex3]{Example 3}: Fracture network geometry, well locations and spatial mesh (left) and close-up of fractures with mesh and well cells (right). The blue fracture cell in \domain[2] marks injection, whereas the red cell in \domain[3] indicates production. }
\label{fig:ex3_geometry}     
\end{figure}

We consider two immersed fractures in a cube shaped domain of side length \SI{1000}{\metre} centred \SI{1500}{\metre} below the surface.
Each fracture contains one well, implemented as a source or sink term in a single cell, with injection in the leftmost fracture, \domain[2], and production in the rightmost fracture, \domain[3]. The domain, fracture geometry and spatial mesh is shown in Fig.~\ref{fig:ex3_geometry}.

The flow rate is \SI{5}{\liter\per\second} for both wells and the injection temperature is \SI{30}{\kelvin} below the formation temperature.
The anisotropic boundary tractions are based on lithostatic stress, with
\begin{gather*}\label{eq:ex3_stress_tensor}
\begin{aligned}
\stress[xx]=0.6\density[s]\gravity z, \qquad \stress[yy]= 1.2 \density[s] \gravity z, \qquad \stress[zz]=\density[s]\gravity z.
\end{aligned}
\end{gather*}
This background stress implies that \domain[3], with normal vector $\normal[3]=[1,0,0]^T$,  is initially more favourably oriented for propagation than is \domain[2].

\begin{figure}
    \includegraphics[width=.98\textwidth]{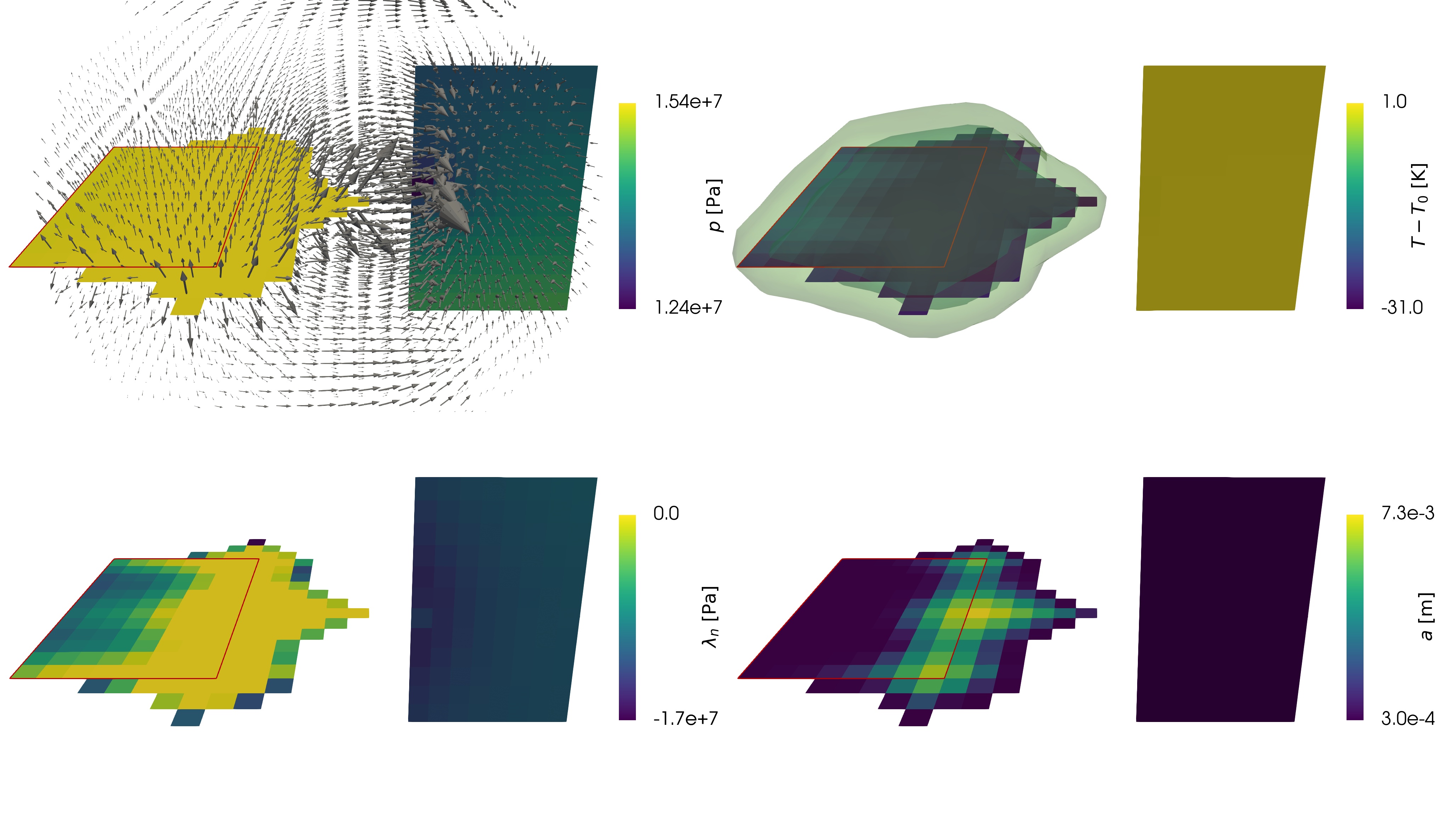}
\caption{\hyperref[sec:ex3]{Example 3}: Solution and fracture geometry at the end of the simulation. \pressure, \aperture and \normalTraction are shown on the fractures, while \temperature is shown both on the fractures and as contour lines indicating where the matrix is significantly cooled (\SIlist{10;20}{\kelvin} below initial formation temperature). The red rectangle in the bottom left figure indicates the initial shape of \domain[2]. 
}
\label{fig:ex3_solutions}     
\end{figure}

\begin{figure}
    \includegraphics[width=.495\textwidth]{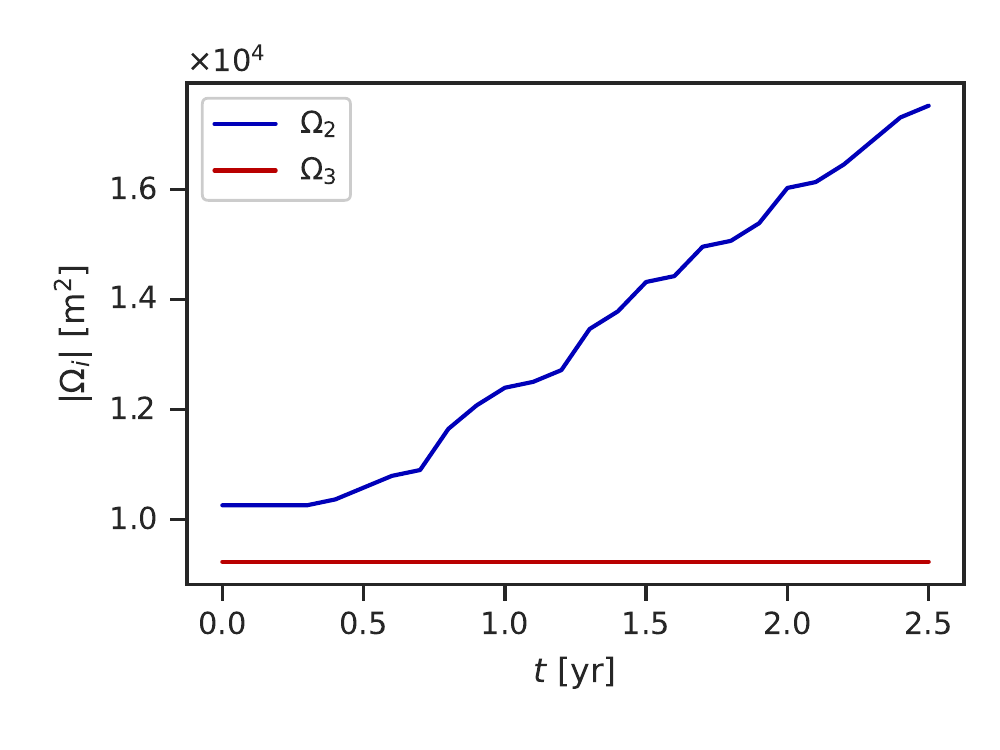}
\includegraphics[width=.495\textwidth]{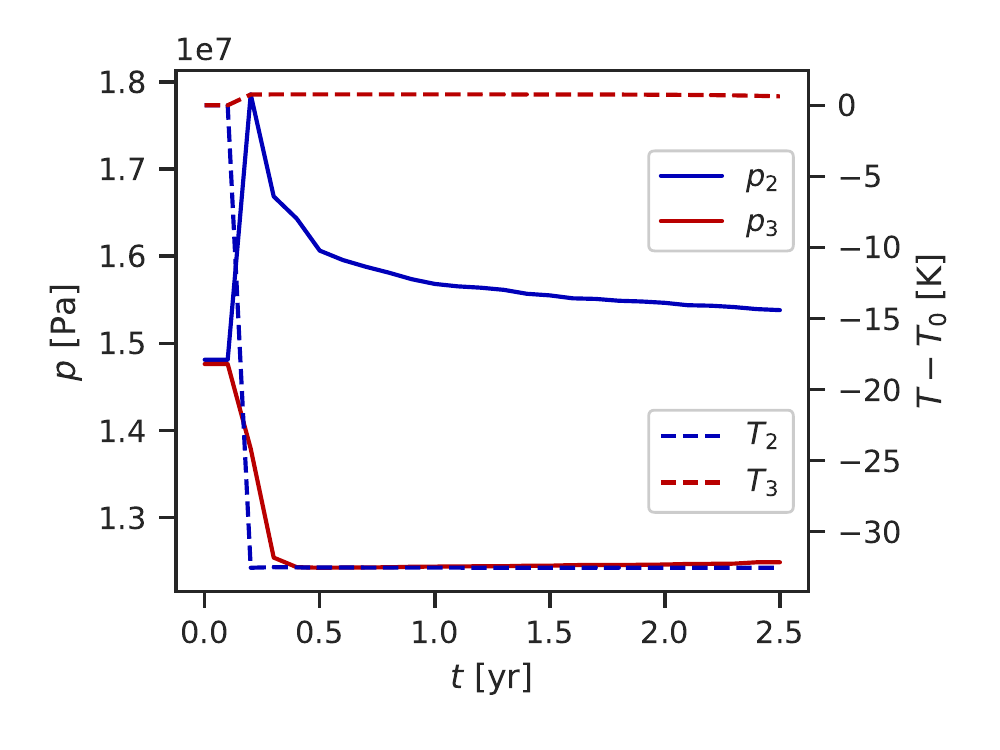}
\caption{\hyperref[sec:ex3]{Example 3}: Size of the two fractures vs. time (left) and pressure and temperature vs. time in the injection well cells (right).
 Only the injection fracture \domain[2] grows.}
\label{fig:ex3_plots}     
\end{figure}

The Fig.~\ref{fig:ex3_plots} fracture size plot shows that \domain[2] grows at a steady speed after an initial phase of limited propagation.
The growth is driven by elevated pressure due to injection and matrix cooling, which is most pronounced on the side of \domain[2] facing \domain[3] due to the advective component of the heat flow, cf.\ Fig.~\ref{fig:ex3_solutions}.
Assuming the thermal driving force to dominate, which is reasonable given the relative size of injection pressure and background stresses, the relatively constant speed could be linked to the constant rate and temperature of injection.

The fracture \domain[3], where fluid is produced, does not propagate at all. Towards the end of the \SI{2.5}{\year} simulation, the magnitude of normal traction on \domain[3] has increased considerably  relative to the initial value of approx. $\SI{1e7}{\pascal}$, see Fig.~\ref{fig:ex3_plots}.
We attribute this  to the contraction ensuing from matrix cooling surrounding \domain[2], which leads to a larger proportion of the compressive forces being supported by the non-cooled surroundings, including \domain[3]. 


Figure \ref{fig:ex3_plots} also shows temperature and pressure in the two wells throughout the simulation. 
Most notably, injection pressure gradually declines.
This increased injectivity in \domain[2] is caused by the combination of an increased aperture in the pre-existing part of the fracture, and the increase in the geometric extension of the fracture.
Thus, fracture deformation caused by thermal and hydraulic stimulation strongly affects the (flow) properties, providing a clear example of the two-way process-structure interaction characteristic of fractured porous media. 

This simulation indicates that long-term cooling during  geothermal energy production may alter the stress state to a stage where fractures  propagate.
It is thus important to develop simulation tools that can incorporate such changes to fracture geometry, in addition to handling multiphysics processes in the reservoir.
Moreover, the injection pressure evolution shows the importance of also capturing deformation of existing fractures in the same model.



\subsubsection{Example 4: Thermal fracturing and natural convection}\label{sec:ex4}
As a final example, we consider fracture propagation driven by cooling that is mainly caused by convection cells inside vertical fractures. 
The process, known as convective downward migration, has been proposed as a mechanism for transport of heat in the deep roots of volcanic geothermal systems \cite{lister1974penetration,bodvarsson1982terrestrial}. It is also predicted to have an important role in the source mechanism of hydrothermal activity in a more general perspective \cite{bodvarsson1982terrestrial,axelsson_1985}.

\begin{figure}
    \includegraphics[width=.48\textwidth]{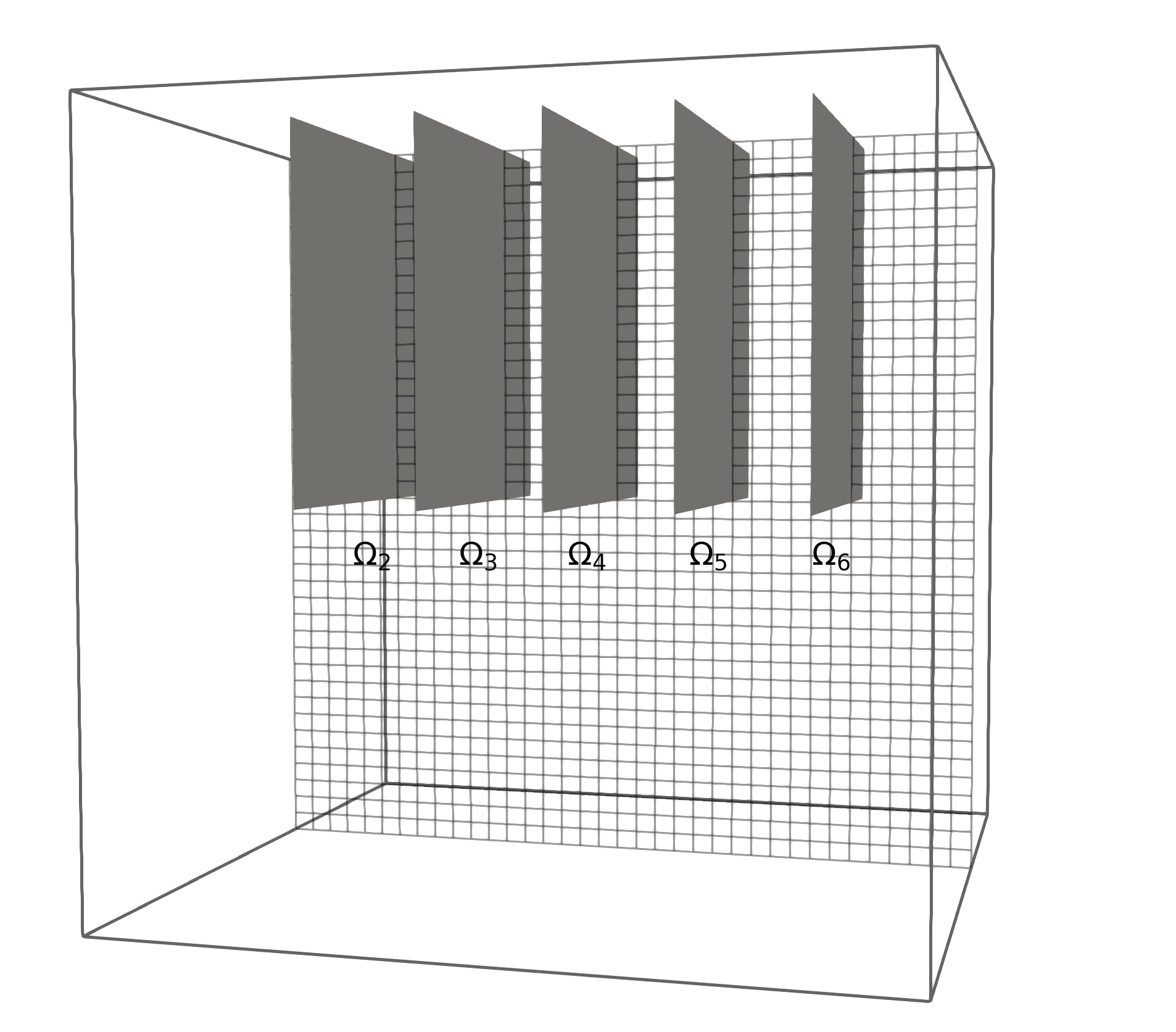}
    \includegraphics[width=.495\textwidth]{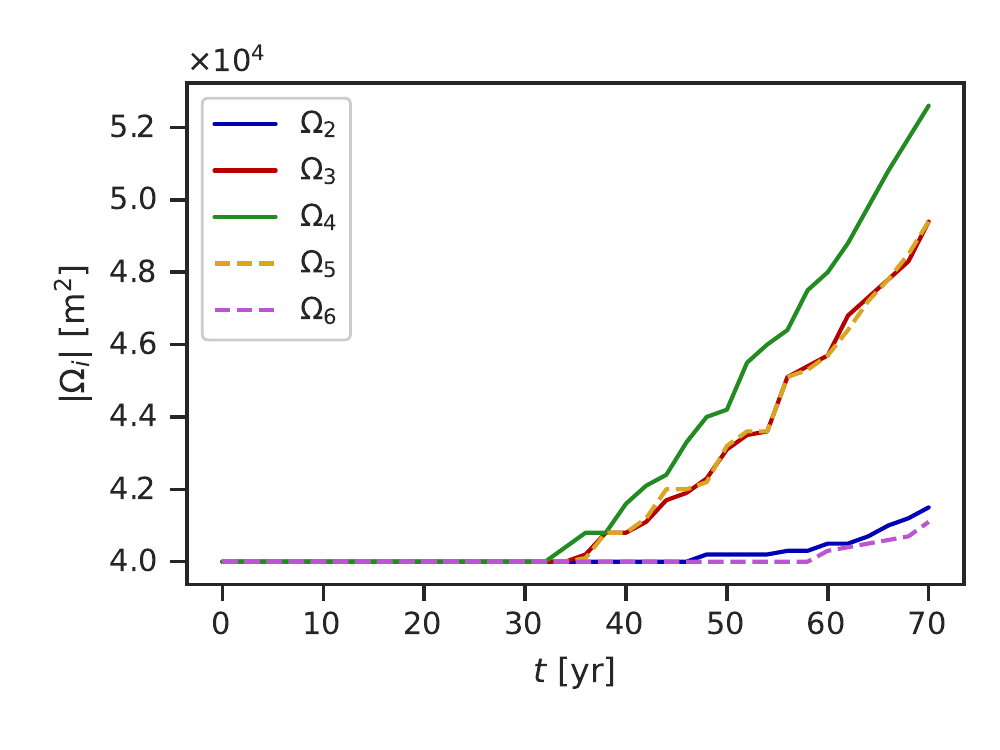}
\caption{\hyperref[sec:ex4]{Example 4}: Initial geometry (left) and size of the five fractures vs. time (right).
}
\label{fig:ex4_geometry_and_sizes}     
\end{figure}

\begin{figure}
    \includegraphics[width=.98\textwidth]{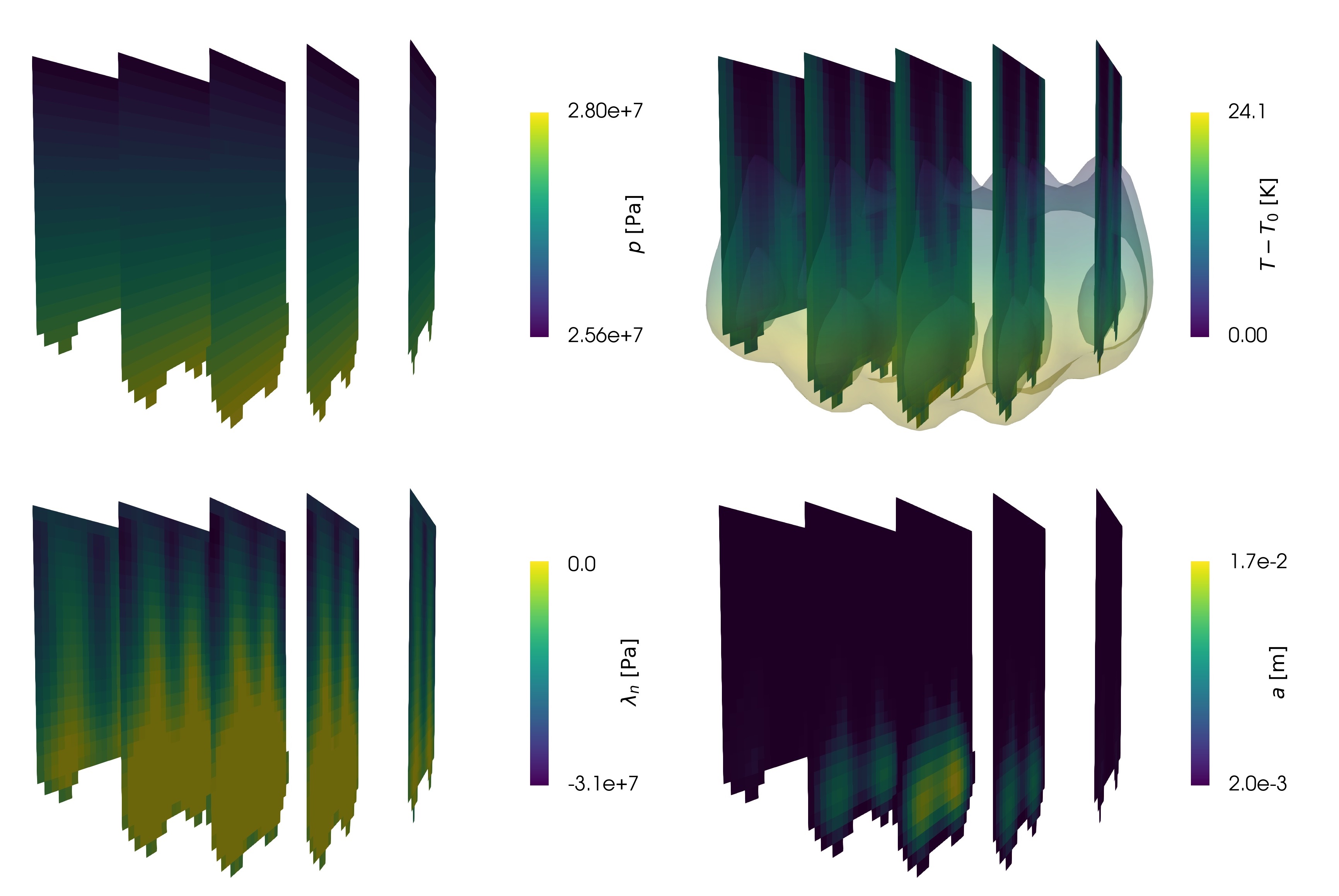}
\caption{\hyperref[sec:ex4]{Example 4}: Solution and fracture geometry at the end of the simulation. \pressure, \aperture and \normalTraction are shown on the fractures, while \temperature is shown both on the fractures and as contour surfaces indicating where the matrix is significantly cooled (\SIlist{7.5;15}{\kelvin} below initial formation temperature). 
}
\label{fig:ex4_solutions}     
\end{figure}

\begin{figure}
     \includegraphics[width=.98\textwidth]{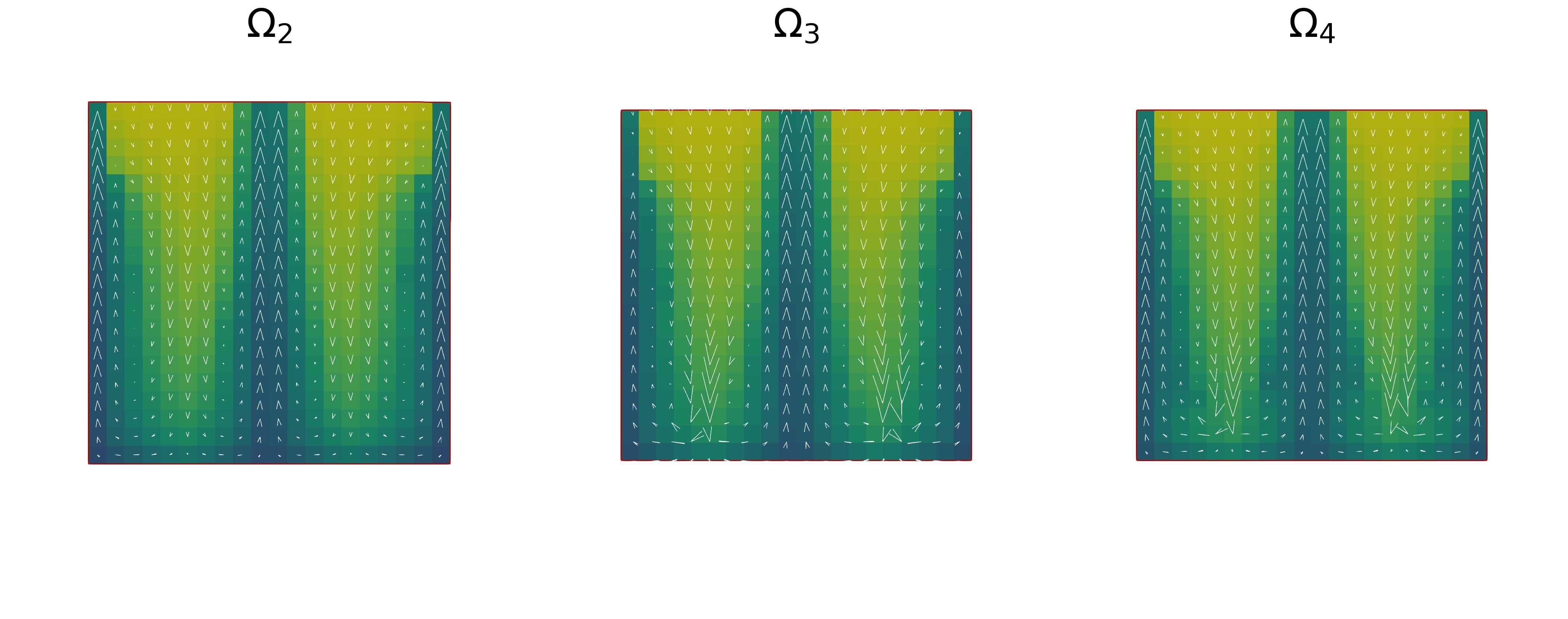}
     \includegraphics[width=.98\textwidth]{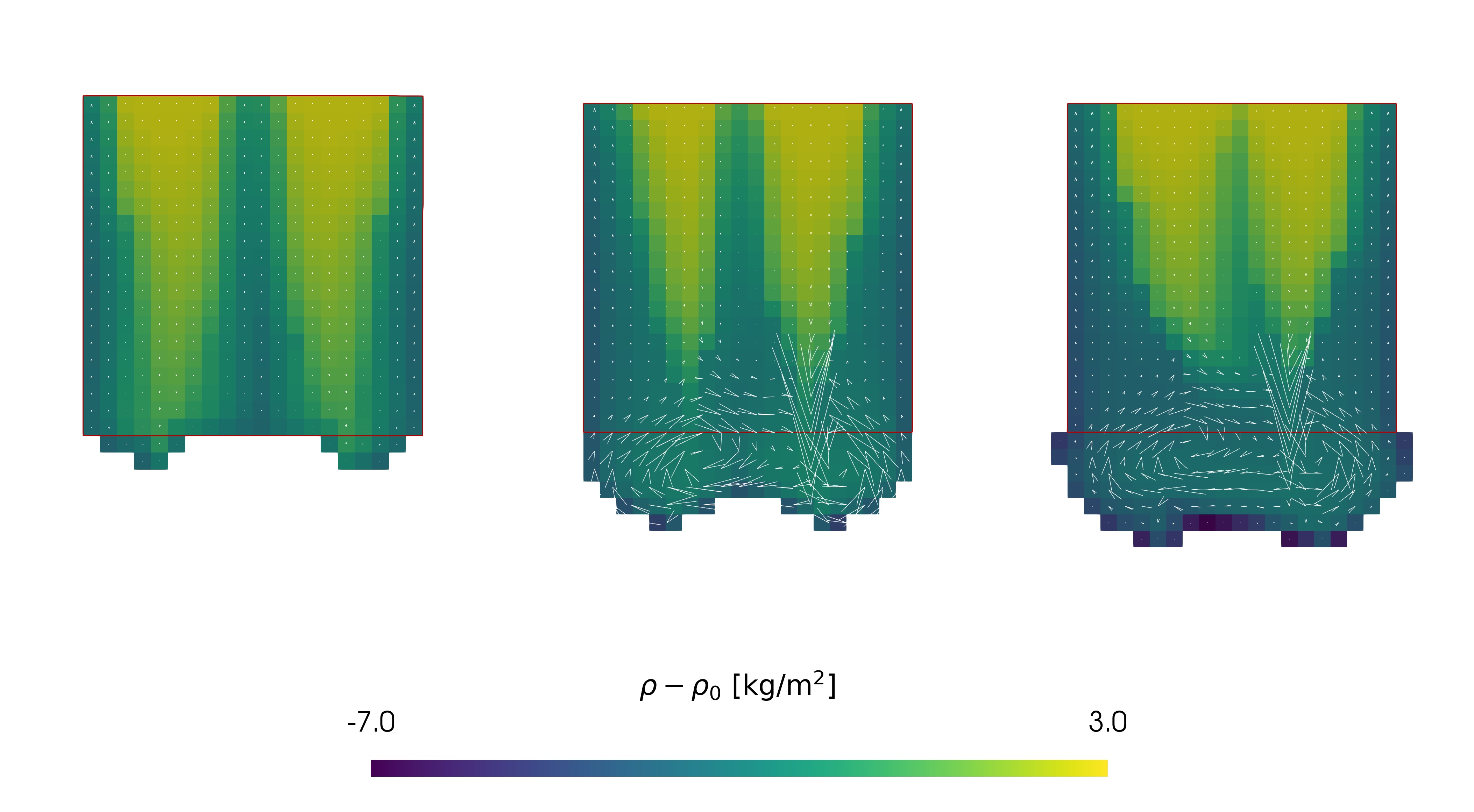}
\caption{\hyperref[sec:ex4]{Example 4}: Flux fields and density distribution for three fractures. The top and bottom row correspond to just before propagation onset ($t=\SI{32}{\year}$)  and  the end of simulation ($t=\SI{70}{\year}$), respectively. The scaling of the arrows indicating flux direction and magnitude is a factor five larger for the top row, when fluxes are smaller due to smaller apertures. The red rectangles indicate initial fracture geometry.}
\label{fig:ex4_fluxes}     
\end{figure}

We consider five vertical fractures evenly spaced along the $x$ direction and extending from the top boundary half-way through the cube-shaped domain with side length \SI{400}{\metre}, see Fig.~\ref{fig:ex4_geometry_and_sizes}. The domain is centred \SI{2800}{\metre} below the surface to mimic conditions within the earth crust where natural heat convection is likely to take place. Boundary and initial conditions are hydrostatic pressure and temperature according to a vertical gradient of \SI{-0.15}{\kelvin\per\meter} and upper boundary temperature \SI{500}{\kelvin}, considered to represent background temperature gradient close to the boundaries of geothermal areas. This value is estimated between \num{-0.10} and \SI{-0.15}{\kelvin\per\meter} within Iceland's active zone of volcanism and rifting, where many high temperature systems exist \cite{agustsson2004TheTO}. The boundary traction is the same as in the previous example and the simulation time is \num{70} years. The results are displayed in Figures \ref{fig:ex4_geometry_and_sizes}, \ref{fig:ex4_solutions} and \ref{fig:ex4_fluxes}.

The vertical temperature gradient leads to instabilities in fluid density, which triggers convection cells inside the fracture, see Fig.~\ref{fig:ex4_fluxes}.
As shown by the  temperature contour surfaces in Fig.~\ref{fig:ex4_solutions}, the resulting energy transport cools the rock surrounding the fractures, to the point where propagation occurs at the lower end of the fractures.
This change in fracture geometry, together with changes in aperture in the existing fracture due to contraction of the surrounding rock, again gives feedback to the fluid convection, as is evident from the difference in flow patterns between the solutions at the two different times reported in Fig.~\ref{fig:ex4_fluxes}. 
As in \hyperref[sec:ex3]{Example 3}, we see evidence of tight process-structure interaction, with the convection-induced cooling altering thermo-poro-mechanical stress sufficiently for the fractures to open and propagate.

Figure~\ref{fig:ex4_geometry_and_sizes} displays size evolution for individual fractures. 
Propagation begins approximately half-way through the simulation, first for the fracture in the center of the domain. 
Even after all fractures have started propagating, the central fractures \domain[3], \domain[4] and \domain[5] propagate significantly faster than the two outermost.
This should be understood in the context of the compressive boundary conditions:
The normal tractions on \domain[2] and \domain[6], respectively, are not relieved by the cooling of any fractures lying between them and the left and right boundary.

After onset, propagation continues  until the end of the simulation, but not at all time steps for all  propagating fractures, and certainly not along the entire propagation front.
This is because the matrix surrounding the new part of a fracture must be cooled before the fracture proceeds, and indicates that the fracture growth is stable as in \hyperref[sec:ex2]{Example 2} and that the propagation speed is resolved in the temporal discretisation.
An approximate downward propagation speed for fractures 3-5 is obtained by dividing the estimated slopes from Fig.~\ref{fig:ex4_geometry_and_sizes} by the initial lateral fracture length \SI{200}{\metre}, yielding $\sim\SI{2}{\metre \per \year}$.
 The setup for this test case is based on average properties in high temperature settings, and the results are in agreement with previous assessments of \SIrange{0.3}{5}{\metre \per \year} \cite{bodvarsson1982terrestrial,bjornsson_penetration_1982}, using a simple relation between the temperature difference sufficient for thermal stress to outweigh the hydrostatic force to keep the fracture closed, at approximatly 3 km depth in the crust with average properties of water and rock similar to the example.  
 



\section{Conclusion}\label{sec:conclusion}
 While both numerical models considering flow and heat transfer in fractured media and models considering deformation of poroelastic media and fracture mechanics have separately been studied extensively, models which combine these fields are more recent.
In the current work, we present a novel numerical model that couples fracture contact mechanics and propagation with deformation, flow and heat transfer in fractured thermo-poroelastic media. The methodology is built on a multi-point control-volume framework, combined with an active-set approach for fracture contact mechanics. The fracture propagation is based on stress intensity factors, and computed using a variant of the displacements correlation method. In the numerical model, fractures are restricted to propagate conforming to the existing grid. The numerical results show mesh convergence for computation of stress-intensity factors and fracture propagation speeds. Focusing on tensile fracture propagation, three-dimensional numerical test cases also show how the model can be used to investigate fracture propagation caused by forced and natural convection, exemplified by long-term thermal reservoir stimulation due to cooling and convective downward migration of fractures. The simulations demonstrate the need for coupled models accounting for both contact mechanics and fracture propagation as well as the coupled thermo-poroelasticity.  

\paragraph{Acknowledgements}
Funding: This work was supported by the Research Council of Norway and Equinor ASA through grants number 267908 and 308733.

\bibliographystyle{spmpsci}      
\bibliography{bibliography}

\begin{thebibliography}{10}
\providecommand{\url}[1]{{#1}}
\providecommand{\urlprefix}{URL }
\expandafter\ifx\csname urlstyle\endcsname\relax
  \providecommand{\doi}[1]{DOI~\discretionary{}{}{}#1}\else
  \providecommand{\doi}{DOI~\discretionary{}{}{}\begingroup
  \urlstyle{rm}\Url}\fi

\bibitem{keilegavlen2020porepyv130}
Porepy v1.3.0 source code.
\newblock DOI
  \href{https://doi.org/10.5281/zenodo.4314343}{10.5281/zenodo.4314343}

\bibitem{Aavatsmark2002}
Aavatsmark, I.: An introduction to multipoint flux approximations for
  quadrilateral grids.
\newblock Comput Geosci \textbf{6}, 405--432 (2002)

\bibitem{agustsson2004TheTO}
{\'A}g{\'u}stsson, K., Fl{\'o}venz, {\'O}.G.: The thickness of the seismogenic
  crust in iceland and its implications for geothermal systems (2005).
\newblock Proceedings World Geothermal Congress 2005, Antalya, Turkey, 24-29
  April 2005

\bibitem{axelsson_1985}
Axelsson, G.: Hydrology and thermomechanics of liquid-dominated hydrothermal
  systems in iceland. (1985).
\newblock PhD-Thesis. Oregon State University, USA, (1985)

\bibitem{baroth_uncertainty_2019}
Baroth, J.: Uncertainty propagation through {Thermo}-{Hydro}-{Mechanical}
  modelling of concrete cracking and leakage – {Application} to containment
  buildings.
\newblock In: Proceedings of the 10th {International} {Conference} on
  {Fracture} {Mechanics} of {Concrete} and {Concrete} {Structures}. IA-FraMCoS
  (2019)

\bibitem{barton1976shear}
Barton, N.: The shear strength of rock and rock joints.
\newblock In: Int J Rock Mech Min Sci Geomech Abstr, vol.~13, pp. 255--279
  (1976)

\bibitem{berge2020hmdfm}
Berge, R.L., Berre, I., Keilegavlen, E., Nordbotten, J.M., Wohlmuth, B.: Finite
  volume discretization for poroelastic media with fractures modeled by contact
  mechanics.
\newblock Int J Numer Methods Eng \textbf{121}(4), 644--663 (2020)

\bibitem{berre2020verification}
Berre, I., Boon, W.M., Flemisch, B., Fumagalli, A., Gl{\"a}ser, D.,
  Keilegavlen, E., Scotti, A., Stefansson, I., Tatomir, A., Brenner, K.,
  et~al.: Verification benchmarks for single-phase flow in three-dimensional
  fractured porous media.
\newblock Adv Water Resour \textbf{147}, 103,759 (2020)

\bibitem{berre2019flow}
Berre, I., Doster, F., Keilegavlen, E.: Flow in fractured porous media: A
  review of conceptual models and discretization approaches.
\newblock Transp Porous Med \textbf{130}(1), 215--236 (2019)

\bibitem{biot1941general}
Biot, M.A.: General theory of three-dimensional consolidation.
\newblock J Appl Phys \textbf{12}(2), 155--164 (1941)

\bibitem{bjornsson_penetration_1982}
Bj{\"o}rnsson, H., Bj{\"o}rnsson, S., Sigurgeirsson, T.: Penetration of water
  into hot rock boundaries of magma at {Gr{\'i}msv{\"o}tn}.
\newblock Nature \textbf{295}(5850), 580--581 (1982)

\bibitem{bjornsson_heatmass1987}
Bj{\"o}rnsson, S., Stef{\'a}nsson, V.: Heat and mass transport in geothermal
  reservoirs.
\newblock In: J.~Bear, M.Y. Corapcioglu (eds.) Advances in Transport Phenomena
  in Porous Media, NATO ASI Series (Series E: Applied Sciences), vol. 128, pp.
  143--183 (1987)

\bibitem{bodvarsson1982terrestrial}
Bodvarsson, G.: Terrestrial energy currents and transfer in iceland.
\newblock Continental and oceanic rifts \textbf{8}, 271--282 (1982)

\bibitem{bois2012use}
Bois, A.P., Garnier, A., Galdiolo, G., Laudet, J.B., et~al.: Use of a
  mechanistic model to forecast cement-sheath integrity.
\newblock SPE Drill Completion \textbf{27}(02), 303--314 (2012)

\bibitem{Boon2018}
Boon, W.M., Nordbotten, J.M., Yotov, I.: Robust discretization of flow in
  fractured porous media.
\newblock SIAM J Numer Anal \textbf{56}(4), 2203--2233 (2018)

\bibitem{de2016gradient}
de~Borst, R., Verhoosel, C.V.: Gradient damage vs phase-field approaches for
  fracture: Similarities and differences.
\newblock Comput Method Appl M \textbf{312}, 78--94 (2016)

\bibitem{both2019gradient}
Both, J.W., Kumar, K., Nordbotten, J.M., Radu, F.A.: The gradient flow
  structures of thermo-poro-visco-elastic processes in porous media.
\newblock arXiv preprint arXiv:1907.03134  (2019)

\bibitem{bouhjiti_statistical_2018}
Bouhjiti, D.E.M., Baroth, J., Briffaut, M., Dufour, F., Masson, B.: Statistical
  modeling of cracking in large concrete structures under
  {Thermo}-{Hydro}-{Mechanical} loads: {Application} to {Nuclear} {Containment}
  {Buildings}. {Part} 1: {Random} field effects (reference analysis).
\newblock Nucl Eng Des \textbf{333}, 196--223 (2018)

\bibitem{bouhjiti_statistical_2018-1}
Bouhjiti, D.E.M., Blasone, M.C., Baroth, J., Dufour, F., Masson, B.,
  Michel-Ponnelle, S.: Statistical modelling of cracking in large concrete
  structures under {Thermo}-{Hydro}-{Mechanical} loads: {Application} to
  {Nuclear} {Containment} {Buildings}. {Part} 2: {Sensitivity} analysis.
\newblock Nucl Eng Des \textbf{334}, 1--23 (2018)

\bibitem{brun2018upscaling}
Brun, M.K., Berre, I., Nordbotten, J.M., Radu, F.A.: Upscaling of the coupling
  of hydromechanical and thermal processes in a quasi-static poroelastic
  medium.
\newblock Transp Porous Med \textbf{124}(1), 137--158 (2018)

\bibitem{chan1970displacementcorrelation}
Chan, S., Tuba, I., Wilson, W.: On the finite element method in linear fracture
  mechanics.
\newblock Eng Frac Mech \textbf{2}(1), 1 -- 17 (1970)

\bibitem{cheng2016poroelasticity}
Cheng, A.H.D.: Poroelasticity, vol.~27.
\newblock Springer (2016)

\bibitem{coussy2004poromechanics}
Coussy, O.: Poromechanics.
\newblock Wiley (2004)

\bibitem{cusini2020simulation}
Cusini, M., White, J.A., Castelletto, N., Settgast, R.R.: Simulation of coupled
  multiphase flow and geomechanics in porous media with embedded discrete
  fractures.
\newblock arXiv preprint arXiv:2007.05069  (2020)

\bibitem{Hau2020}
Dang-Trung, H., Keilegavlen, E., Berre, I.: Numerical modeling of wing crack
  propagation accounting for fracture contact mechanics.
\newblock Int J Solids Struct \textbf{204-205}, 233 -- 247 (2020)

\bibitem{umfpack}
Davis, T.A.: Algorithm 832: Umfpack v4.3---an unsymmetric-pattern multifrontal
  method.
\newblock ACM Trans. Math. Softw. \textbf{30}(2), 196–199 (2004)

\bibitem{deb2020extended}
Deb, R., Jenny, P.: An extended finite volume method and fixed-stress approach
  for modeling fluid injection--induced tensile opening in fractured
  reservoirs.
\newblock Int J Numer Anal Methods Geomech \textbf{44}(8), 1128--1144 (2020)

\bibitem{flemisch2018benchmarks}
Flemisch, B., Berre, I., Boon, W., Fumagalli, A., Schwenck, N., Scotti, A.,
  Stefansson, I., Tatomir, A.: Benchmarks for single-phase flow in fractured
  porous media.
\newblock Adv. Water Resour. \textbf{111}, 239--258 (2018)

\bibitem{franceschini2019block}
Franceschini, A., Castelletto, N., Ferronato, M.: Block preconditioning for
  fault/fracture mechanics saddle-point problems.
\newblock Comput Method Appl M \textbf{344}, 376--401 (2019)

\bibitem{franceschini2020algebraically}
Franceschini, A., Castelletto, N., White, J.A., Tchelepi, H.A.: Algebraically
  stabilized lagrange multiplier method for frictional contact mechanics with
  hydraulically active fractures.
\newblock Comput Method Appl M \textbf{368}, 113,161 (2020)

\bibitem{gallyamov2018discrete}
Gallyamov, E., Garipov, T., Voskov, D., Van~den Hoek, P.: Discrete fracture
  model for simulating waterflooding processes under fracturing conditions.
\newblock Int J Numer Anal Methods Geomech \textbf{42}(13), 1445--1470 (2018)

\bibitem{gao_three-dimensional_2020}
Gao, Q., Ghassemi, A.: Three-{Dimensional} {Thermo}-{Poroelastic} {Modeling}
  and {Analysis} of {Flow}, {Heat} {Transport} and {Deformation} in {Fractured}
  {Rock} with {Applications} to a {Lab}-{Scale} {Geothermal} {System}.
\newblock Rock Mech Rock Eng \textbf{53}(4), 1565--1586 (2020)

\bibitem{garipov2019discrete}
Garipov, T., Hui, M.: Discrete fracture modeling approach for simulating
  coupled thermo-hydro-mechanical effects in fractured reservoirs.
\newblock Int J Rock Mech Min \textbf{122}, 104,075 (2019)

\bibitem{garipov2016discrete}
Garipov, T., Karimi-Fard, M., Tchelepi, H.: Discrete fracture model for coupled
  flow and geomechanics.
\newblock Comput Geosci \textbf{20}(1), 149--160 (2016)

\bibitem{giovanardi2017unfitted}
Giovanardi, B., Formaggia, L., Scotti, A., Zunino, P.: Unfitted fem for
  modelling the interaction of multiple fractures in a poroelastic medium.
\newblock In: Geometrically Unfitted Finite Element Methods and Applications,
  pp. 331--352. Springer (2017)

\bibitem{griffith1921rupture}
Griffith, A.A.: Vi. the phenomena of rupture and flow in solids.
\newblock Phil T R Soc. A \textbf{221}(582-593), 163--198 (1921)

\bibitem{hueber2008paper}
H{\"u}eber, S., Stadler, G., Wohlmuth, B.I.: A primal-dual active set algorithm
  for three-dimensional contact problems with coulomb friction.
\newblock SIAM J Sci Comput \textbf{30}(2), 572--596 (2008)

\bibitem{jaffre2011interdimupwind}
Jaffré, J., Mnejja, M., Roberts, J.: A discrete fracture model for two-phase
  flow with matrix-fracture interaction.
\newblock Procedia Computer Science \textbf{4}, 967 -- 973 (2011).
\newblock Proceedings of the International Conference on Computational Science,
  ICCS 2011

\bibitem{Karimi-Fard2003efficient}
Karimi-Fard, M., Durlofsky, L.J., Aziz, K.: An {E}fficient
  {D}iscrete-{F}racture {M}odel {A}pplicable for {G}eneral-{P}urpose
  {R}eservoir {S}imulators.
\newblock SPE J \textbf{9}(2), 227--236 (2003)

\bibitem{keilegavlen2020porepy}
Keilegavlen, E., Berge, R., Fumagalli, A., Starnoni, M., Stefansson, I.,
  Varela, J., Berre, I.: Porepy: An open-source software for simulation of
  multiphysics processes in fractured porous media.
\newblock Comput Geosci  (2020).
\newblock \doi{10.1007/s10596-020-10002-5}

\bibitem{keilegavlen2017finite}
Keilegavlen, E., Nordbotten, J.M.: Finite volume methods for elasticity with
  weak symmetry.
\newblock Int J Numer Methods Eng \textbf{112}(8), 939--962 (2017)

\bibitem{khoei2014mesh}
Khoei, A.R., Vahab, M., Haghighat, E., Moallemi, S.: A mesh-independent finite
  element formulation for modeling crack growth in saturated porous media based
  on an enriched-fem technique.
\newblock Int J Fracture \textbf{188}(1), 79--108 (2014)

\bibitem{kiraly1988mixedfem}
Kir{\'a}ly, L.: Large scale 3-d groundwater flow modelling in highly
  heterogeneous geologic medium.
\newblock In: Groundwater flow and quality modelling, pp. 761--775. Springer
  (1988)

\bibitem{kogbara2013review}
Kogbara, R.B., Iyengar, S.R., Grasley, Z.C., Masad, E.A., Zollinger, D.G.: A
  review of concrete properties at cryogenic temperatures: Towards direct lng
  containment.
\newblock Constr Build Mater \textbf{47}, 760--770 (2013)

\bibitem{lin_integrity_2020}
Lin, Y., Deng, K., Yi, H., Zeng, D., Tang, L., Wei, Q.: Integrity tests of
  cement sheath for shale gas wells under strong alternating thermal loads.
\newblock Natural Gas Industry B  (2020)

\bibitem{lister1974penetration}
Lister, C.: On the penetration of water into hot rock.
\newblock Geophys J Int \textbf{39}(3), 465--509 (1974)

\bibitem{Martin2005}
Martin, V., Jaffr{\'e}, J., Roberts, J.E.: Modeling {F}ractures and {B}arriers
  as {I}nterfaces for {F}low in {P}orous {M}edia.
\newblock SIAM J Sci Comput \textbf{26}(5), 1667--1691 (2005)

\bibitem{nejati2015}
Nejati, M., Paluszny, A., Zimmerman, R.W.: On the use of quarter-point
  tetrahedral finite elements in linear elastic fracture mechanics.
\newblock Eng Frac Mech \textbf{144}, 194 -- 221 (2015)

\bibitem{nordbotten2016biot}
Nordbotten, J.: Stable cell-centered finite volume discretization for biot
  equations.
\newblock SIAM J Numer Anal \textbf{54}, 942--968 (2016)

\bibitem{nordbotten2020mpxa}
Nordbotten, J., Keilegavlen, E.: An introduction to multi-point flux (mpfa) and
  stress (mpsa) finite volume methods for thermo-poroelasticity.
\newblock arXiv preprint arXiv:2001.01990  (2020)

\bibitem{paluszny2011}
Paluszny, A., Zimmerman, R.W.: Numerical simulation of multiple 3d fracture
  propagation using arbitrary meshes.
\newblock Comput Method Appl M \textbf{200}(9), 953 -- 966 (2011)

\bibitem{ren2016fully}
Ren, G., Jiang, J., Younis, R.M.: A fully coupled xfem-edfm model for
  multiphase flow and geomechanics in fractured tight gas reservoirs.
\newblock Procedia Computer Science \textbf{80}, 1404--1415 (2016)

\bibitem{richard2005fracturepropagation}
Richard, H.A., Fulland, M., Sander, M.: Theoretical crack path prediction.
\newblock Fatigue Fract Eng M \textbf{28}(1-2), 3--12 (2005)

\bibitem{salimzadeh2018thm}
Salimzadeh, S., Paluszny, A., Nick, H.M., Zimmerman, R.W.: A three-dimensional
  coupled thermo-hydro-mechanical model for deformable fractured geothermal
  systems.
\newblock Geothermics \textbf{71}, 212 -- 224 (2018)

\bibitem{salimzadeh2017three}
Salimzadeh, S., Paluszny, A., Zimmerman, R.W.: Three-dimensional poroelastic
  effects during hydraulic fracturing in permeable rocks.
\newblock Int J Solids Struct \textbf{108}, 153--163 (2017)

\bibitem{selvadurai2017thermo}
Selvadurai, A.P., Suvorov, A.: Thermo-poroelasticity and geomechanics.
\newblock Cambridge University Press (2017)

\bibitem{settgast2017fully}
Settgast, R.R., Fu, P., Walsh, S.D., White, J.A., Annavarapu, C., Ryerson,
  F.J.: A fully coupled method for massively parallel simulation of
  hydraulically driven fractures in 3-dimensions.
\newblock Int J Numer Anal Method Geomech \textbf{41}(5), 627--653 (2017)

\bibitem{siratovich2015saturated}
Siratovich, P.A., Villeneuve, M.C., Cole, J.W., Kennedy, B.M., B{\'e}gu{\'e},
  F.: Saturated heating and quenching of three crustal rocks and implications
  for thermal stimulation of permeability in geothermal reservoirs.
\newblock Int J Rock Mech Min \textbf{80}, 265--280 (2015)

\bibitem{sneddon1946analytical}
Sneddon, I.N.: The distribution of stress in the neighbourhood of a crack in an
  elastic solid.
\newblock Philos T R Soc A \textbf{187}(1009), 229--260 (1946)

\bibitem{stefansson2020thm}
Stefansson, I., Berre, I., Keilegavlen, E.: A fully coupled numerical model of
  thermo-hydro-mechanical processes and fracture contact mechanics in porous
  media.
\newblock arXiv preprint arXiv:2008.06289  (2020)

\bibitem{stefansson2020sourcecode_propagation}
Stefansson, I., Keilegavlen, E.: Run scripts for {PorePy} simulations.
\newblock DOI
  \href{https://doi.org/10.5281/zenodo.4316328}{10.5281/zenodo.4316328} (2020)

\bibitem{ucar2018finite}
Ucar, E., Keilegavlen, E., Berre, I., Nordbotten, J.M.: A finite-volume
  discretization for deformation of fractured media.
\newblock Comput Geosci \textbf{22}(4), 993--1007 (2018)

\bibitem{wohlmuth2011}
Wohlmuth, B.: Variationally consistent discretization schemes and numerical
  algorithms for contact problems.
\newblock Acta Numer \textbf{20}, 569--734 (2011)

\bibitem{wu_investigation_2020}
Wu, Z., Zhou, Y., Weng, L., Liu, Q., Xiao, Y.: Investigation of thermal-induced
  damage in fractured rock mass by coupled {FEM}-{DEM} method.
\newblock Comput Geosci  (2020)

\bibitem{zimmerman1996cubic}
Zimmerman, R.W., Bodvarsson, G.S.: Hydraulic conductivity of rock fractures.
\newblock Transp Porous Med \textbf{23}(1), 1--30 (1996)

\end{thebibliography}
\end{document}